\documentclass[12pt]{article}
\usepackage[square,sort,comma,numbers]{natbib}
\usepackage[utf8]{inputenc}
\usepackage[english]{babel}
\usepackage{color}
\usepackage[usenames,dvipsnames,svgnames,table]{xcolor}
\usepackage[colorlinks=true, linkcolor=blue, citecolor= red, urlcolor=red]{hyperref}
\usepackage{amssymb,amsthm}
\newtheorem{thm}{Theorem}
\newtheorem*{remark}{Remark}
\usepackage{mathtools}
\usepackage[T1]{fontenc}
\usepackage{amsmath}
\DeclareMathOperator{\tr}{tr}
\usepackage{enumitem}

\usepackage{mathtools}
\usepackage{caption}
\usepackage{layout}
\usepackage{float}
\usepackage{graphicx}
\usepackage[margin = 1in]{geometry}
\setcitestyle{square}

\usepackage{authblk}

\title{Density Power Downweighting and Robust Inference: Some New Strategies}
\author[1]{Saptarshi Roy }
\author[2]{Kaustav Chakraborty}
\author[3]{Somnath Bhadra}
\author[4]{Ayanendranath Basu}
\affil[1]{\textit{Department of Statistics, University of Michigan, Ann Arbor, MI 48109, USA;}} 
\affil[2]{\textit{Department of Statistics, University of Illinois, Urbana-Champaign, IL 61820, USA;}} 
\affil[3]{\textit{Department of Statistics, University of Florida, Gainesville, FL 32611, USA;}} 
\affil[4]{\textit{Interdisciplinary Statistical Research Unit, Indian Statistical Institute, 203 B. T. Road, Kolkata 700108, India.}}

\date{}
\begin{document}
\maketitle




\begin{abstract}
    Preserving the robustness of the procedure has, at the present time, become almost a default requirement for statistical data analysis. Since efficiency at the model and robustness under misspecification of the model are often in conflict, it is important to choose such inference procedures which provide the best compromise between these two concepts. Some minimum Bregman divergence estimators and related tests of hypothesis seem to be able to do well in this respect, with the procedures based on the density power divergence providing the existing standard. In this paper we propose a new family of Bregman divergences which is a superfamily encompassing the density power divergence. This paper describes the inference procedures resulting from this new family of divergences, and makes a strong case for the utility of this divergence family in statistical inference.
\end{abstract}





\section{Introduction}
 In statistical modeling, parameter estimation is an inevitable and formidable task. Accurate estimation of the model facilitates the characterization and the subsequent understanding of the mechanism that generates the observed data.~Statistical distances can be useful tools for the estimation of the model parameters.
 
 Statistical distances can be naturally applied to the case of parametric statistical inference.~The most important idea in parametric minimum distance inference is the quantification of the degree of closeness between the sample data and parametric model as a function of an unknown set of parameters through a suitable distance-like measure. Thus the estimate of the parameter is obtained by minimizing this ``distance'' over the parameter space.
 
 It is worthwhile to mention here that the class of distances which we will consider are not mathematical metrics in the strict sense of the term. They may not be symmetric in their arguments and may not satisfy the triangle inequality. The only properties that we require of these measures are that they should be nonnegative, and should equal zero if and only if the arguments are identically equal. However, we will, somewhat loosely, continue to call them distances, or ``statistical distances''. In a practical sense, the word ``divergence'' is a good descriptor of these measures.  We will, in fact, use the ``minimum distance'' and the ``minimum divergence'' terminologies interchangeably.
 
 Density-based divergences form a special class of statistical distances. Several minimum distance estimators in this family have high model efficiency. In particular, the maximum likelihood estimator (MLE) also belongs to the class of density-based minimum distance estimators, being the minimizer of the likelihood disparity (Lindsay, 1994), which is a version of the Kullback-Leibler divergence. But one of the major drawbacks of the MLE is that it is notoriously nonrobust and even a small proportion of outlying observations can lead to meaningless inference. In fact it is the failure of the classical methods like maximum likelihood to deal with outliers and mild deviations from the model which had led to the emergence of the field of robustness; see, for example, \cite{huber1981robust}, \cite{hampel1986robust}, \cite{maronna2006robust} and \cite{Basu}. However, some of the other members of the class of minimum distance estimators have been observed to do much better in the sense of combining strong robustness with high model efficiency. See, for example, \cite{Csiszar}, \cite{Ali}, \cite{lindsay1994efficiency}, \cite{pardo2005statistical} and \cite{Basu} for a description of the $\phi$-divergence class of minimum distance measures.

A more modern class of minimum distance estimators is based on the family of Bregman divergences. The Bregman divergence (Bregman, 1967) is a distance like measure between points and has been used in mathematics and information theory for some time. When the points are represented by probability distributions, the corresponding Bregman divergence is a statistical distance. See, for example, \cite{jones1990general}, \cite{csiszar1991least}, \cite{banerjee2005clustering} and \cite{stummer2012bregman} for some examples of statistical and related applications of the Bregman divergence. The principal representatives of Bregman divergence estimators in the current statistical literature are the minimum density power divergence estimators (MDPDEs), based on the density power divergence (DPD) class of \cite{bhhj}. Over the last two decades, this class of divergences has provided a popular and frequently used method to balance the trade-off between robustness and efficiency in parameter estimation, hypothesis testing, and related inference.~The minimum divergence estimators based on the DPD have been shown to provide a high degree of stability under model misspecification, often with minimal loss in model efficiency. Our primary purpose in this paper is to refine the minimum distance procedure based on the DPD, so as to achieve even better compromise between efficiency and robustness.
\section{The Bregman Divergence}
Consider a  parametric family of densities $\mathcal{F}=\{f_\theta\; :\theta\in\Theta\subseteq \mathbb{R}^p\}$.
Let $X_1,X_2,\ldots ,X_n$ be i.i.d.~observations from a distribution $G$ having probability density function (pdf) $g$. For the sake of a unified notation we will continue to use the term pdf irrespective of whether the distribution of $G$ is continuous or discrete. Let the common support of  $g$ and $f_\theta$ be $\mathcal{X}\subseteq \mathbb{R}$.
The Bregman divergence between the density $g$ and model density $f_\theta
$ is given by
\begin{equation}
\begin{split}
&D_B(g,f_\theta)\\
&=\int_{\mathcal{X}} \big(B(g(x))-B(f_\theta(x))-(g(x)-f_\theta(x))B'(f_\theta(x))\big)dx,\\ 
\end{split}
\label{1}
\end{equation}
where the index function $B(\cdot)$ is strictly convex and $B'(\cdot)$ represents its first derivative with respect to its argument. In practice, where $f_\theta$ is the pdf of the parametric family, $g$ is the true density, the minimization of the above divergence over the parameter space $\Theta$ will generate the corresponding minimum distance functional which can lead to meaningful inference, depending on the form of the function $B(\cdot)$. The DPD, defined later in this section, is a special case of the Bregman divergence for $\displaystyle B(y)=\frac{y^{1+\alpha}}{\alpha},\hspace{0.1cm} \alpha \geq 0$. 

When the model is differentiable, the general estimating equation under the divergence in Eq.~\ref{1} is
\begin{equation}
\begin{split}
\int u_\theta(x) B^{\prime\prime}(f_\theta(x))f^2_\theta(x)dx-\int u_\theta(x) B^{\prime\prime}(f_\theta(x))&f_\theta(x)g(x)dx\\
&=0,
\end{split}
\label{2}
\end{equation}
or equivalently
\begin{equation}
\begin{split}
\int u_\theta(x) B^{\prime\prime}(f_\theta(x))f^2_\theta(x)dx-\int u_\theta(x) B^{\prime\prime}(f_\theta(x))&f_\theta(x)dG(x)\\
&=0,
\end{split}
\label{3}
\end{equation}
where $u_\theta(x)=\nabla \log f_\theta(x)$ is the score function of the model $f_\theta(x)$, $\nabla$ represents derivative with respect to $\theta$ and $B^{\prime\prime}(\cdot)$ represents the second derivative of $B(\cdot)$ with respect to its argument. Since $G$ is unknown, we construct an empirical version of the divergence in Eq.~\ref{1}, or the estimating equation given in Eq.~\ref{3}, by replacing $G$ (the true data generating distribution) by its empirical counterpart $G_n$. This leads to a class of unbiased (under the model) estimating equations   
\begin{equation}
\begin{split}
\int u_\theta(x) &B^{\prime\prime}(f_\theta(x))f^2_\theta(x)dx\\ &-\frac{1}{n}\sum_{i=1}^{n} u_\theta(X_i) B^{\prime\prime}(f_\theta(X_i))f_\theta(X_i)=0.
\end{split}
\label{4}
\end{equation}
The root of the Eq.~\ref{4} is defined to be the minimum Bregman divergence estimator (MBDE). Here the robustness of the corresponding minimum distance estimator may be at least partially understood by observing the effect of the downweighting function $B''(f_\theta(x))f_\theta(x)$ on $u_\theta(x)$ for less probable values of $x$ under $f_\theta$.~For the DPD, this weight becomes $(\alpha+1)f^\alpha_\theta(x)$.

In this paper we attempt to find a refinement of the DPD downweighting scheme, and, by reconstruction, a corresponding divergence.~We will show that the corresponding minimum distance procedure provides a better compromise between robustness and efficiency in many cases compared to the minimum density power divergence estimator (MDPDE).
\subsection{The Density Power Divergence}
As mentioned earlier, the density power divergence (DPD) is obtained by substituting $B(y)=\dfrac{y^{1+\alpha}}{\alpha}$ in Eq.~\ref{1}. The general form of this divergence, as a function of a nonnegative tuning parameter $\alpha$, is 
\begin{equation}
\textsc{DPD}_{\alpha}(g,f_{\theta})=\int {\Big\{f_{\theta}^{1+\alpha}-\Big(1+\frac{1}{\alpha}\Big)gf_{\theta}^{\alpha}+\frac{1}{\alpha}g^{1+\alpha}\Big\}}.
\label{5}
\end{equation}
For simplicity we have dropped the dummy variable in the above equation. One can define the minimum DPD functional $T_\alpha(G)$ at $G$ through the relation
\begin{equation}
\textsc{DPD}_\alpha(g,f_{T_{\alpha}(G)})=\inf_{\theta\in\Theta}\textsc{DPD}_\alpha(g,f_{\theta}).
\label{6}
\end{equation} 
Under the estimation set up of this paper, the empirical objective function, ignoring the terms independent of $\theta$,  becomes
$$\int f_\theta^{1+\alpha} -\left(1+\frac{1}{\alpha}\right)\frac{1}{n} \sum_{i=1}^{n} f_\theta^\alpha(X_i),$$
and under differentiability of the model, the estimating equation becomes (by equating the negative of the derivative of the above objective function to 0)

\begin{equation}
\centering
\frac{1}{n}\sum_{i=1}^{n} u_\theta(X_i) f_\theta^{\alpha}(X_i)-\int u_\theta(x) f_\theta^{1+\alpha}(x)dx=0.
\label{7}
\end{equation}
It is evident that as $\alpha\rightarrow 0^{+}$,  Eq.~\ref{7} converges to 
the maximum likelihood score equation
\begin{equation}
\hspace{0.9in}\frac{1}{n}\sum_{i=1}^{n}u_\theta(X_i)=0.
\label{8}
\end{equation}
Note that in the part involving real data in Eq.~\ref{7}, a downweighting effect is exerted on the score function 
$u_\theta(x)$ by the factor $f_\theta^\alpha(x)$. This downweighting philosophy will be crucial for developing the new class of procedures. Note that there is no downweighting for the case $\alpha = 0$.
 
The asymptotic properties of the MDPDE have been well studied, and are available, for example, in \cite{Basu}, where the asymptotic distribution of the MDPDE has been explicitly derived. It is useful to note that the MDPDE solves an estimating equation of the form $\sum_{i=1}^{n}\psi(X_i,\theta)=0$, where 
\begin{equation}
\psi(x,\theta)=u_\theta(x) f_\theta^\alpha (x)-\int u_\theta(x) f_\theta ^{1+\alpha}(x)dx.
\label{13}
\end{equation}
 Hence it belongs the class of M-estimators.~So, the asymptotic properties of the MDPDE also follow from M-estimator theory.
 
 \section{A New Divergence}
 Our key philosophy for constructing new divergences and estimation strategies involves manipulating the downweighting factor $B^{\prime\prime}(f_\theta(X_i))f_\theta(X_i)$ in Eq.~\ref{4}. Here we are going to develop a stronger downweighting effect compared to the MDPD estimating equation. Our exploration will generate an estimation scheme with two tuning parameters and we will explore the possibility of coming up with specific candidates which might beat the MDPDEs both in terms of efficiency and robustness.  
\subsection{Choosing the B Function}
The downweighting effect on the score $u_\theta(x)$ applied by the MDPD estimating equation is $f_\theta^\alpha(x)$. As we want to impose a stronger downweighting in relation to this, we wish to choose the $B$ function (or rather, the $B''$ function) so that 
as $x \rightarrow 0^+$, $x B''(x)$ converges to zero faster than $\displaystyle x^\alpha$ for $\alpha > 0$ fixed. (Note, from Eq.~\ref{4}, the downweighting term for $u_\theta(x)$ in the general Bregman divergence is $f_\theta B^{\prime\prime}(f_\theta)$). In particular, we will assume the following conditions on $B''$.

\vspace{0.1in}
\noindent (P1) $B^{\prime\prime}(x)>0 \;\forall \;x>0$, so that $B$ is a strictly convex function over $\mathbb{R}^{+}$.

\vspace{0.1in}
\noindent (P2) $xB^{\prime\prime}(x)$ is an increasing function over $x$  in $(0,\infty)$. Thus the less likely observations will be downweighted more.

\vspace{0.1in}
\noindent (P3) For all $\beta \in (0, 1)$, $\displaystyle \lim_{x\to 0^{+}}\frac{xB^{\prime\prime}(x)}{x^\beta}=0$, i.e., the Bregman formulation attaches weights to the score function which go to zero at a rate faster than the corresponding weights in the MDPD estimating equation.

\vspace{0.1in}
\noindent (P4) $B^{\prime\prime}(x)=x^{\beta}\phi(x,\gamma)$ ($0<\beta\leq 1$, $0<\gamma\leq 1$). Where $\phi(x,\gamma)$ is a continuous and positive function over $\gamma\in (0,1]$ and $x>0$. Furthermore, we demand $\lim_{\gamma\rightarrow 0^{+}}\phi(x,\gamma)=\dfrac{1}{x}$.

To prove that such choice of $B(\cdot)$ satisfying (P1)-(P4) can help us generate divergences which have the desired properties and provide superior inference compared to the DPD, let us first demonstrate the general asymptotic properties of the minimum Bregman divergence estimators. For ease of representation, we refer to the divergence generated by the $B(\cdot)$ function satisfying (P1) to (P4) as $\phi$-DPD.

\subsubsection{General Asymptotic Properties of the MBDE}
We need some regularity assumptions to prove the asymptotic properties of the general MBDE, which we list below.

\vspace{0.1in}
\noindent (A1) The pdfs $f_\theta$ of $X$ have common support, so that the set $\mathcal{X}=\{x|f_\theta (x)>0\}$ is independent of $\theta$. The distribution $G$ is also supported on $\mathcal{X}$, on which the corresponding density $g$ is greater than zero.

\vspace{0.1in}
\noindent (A2) There is an open subset $\omega$ of the parameter space $\Theta$, containing the best fitting parameter $\theta^{g}$ ($\textsc{D}_{B}(g,f_{\theta^{g}})=\inf_{\theta\in \Theta} \textsc{D}_{B}(g,f_\theta)$) such that for almost all $x\in \mathcal{X}$, and all $\theta\in \omega$, the density $f_\theta(x)$ is three times differentiable with respect to $\theta$ and the third partial derivatives are continuous with respect to $\theta$. (The best fitting parameter $\theta^g$ depends on the index function $B(\cdot)$ also, but we suppress that notation for brevity).

\vspace{0.1in}
\noindent (A3) The integrals $\int B^{\prime\prime}(f_\theta(x))f_\theta^2(x)dx$ and $\int B^{\prime\prime}(f_\theta(x))f_\theta(x)g(x) dx$ can be differentiated with respect to $\theta$, and the derivatives can be taken under the integral sign.

\vspace{0.1in}
\noindent (A4) The $p\times p$ matrix $J_B(\theta)$ defined by
\begin{equation*}
\begin{split}
    &J_{B,kl}(\theta)=\\
    &E_g\Big\{\nabla_{kl}\Big(\int [B^{\prime}(f_\theta(x))f_{\theta}(x)-B(f_\theta(x))]dx-B^{\prime}(f_\theta(X))\Big)\Big\}
    \end{split}
\end{equation*}
is positive definite where $E_g$ represents the expectation under the density $g$.~Where $\nabla_{kl}$ represents the partial derivative with respect to the indicated components of $\theta$.

\vspace{0.1in}
\noindent (A5) There exists functions $M_{jkl}(x)$, $j,k,l=1,\ldots, p$, such that
\begin{equation*}
   \begin{split} 
 &\Big|\nabla_{jkl}\Big(\int [B^{\prime}(f_\theta(x))f_{\theta}(x)-B(f_\theta(x))]dx-B^{\prime}(f_\theta(X))\Big)\Big|\\
 &\leq M_{jkl}(X);\;\; \forall \theta \in \omega
 \end{split}
 \end{equation*}
where $E_g[M_{jkl}(X)]<m_{jkl}<\infty\; \forall j,\;k,\;l.$

\begin{thm}
Under the conditions (A1)-(A5), the following results hold.\\
(a) The MBDE estimating equation given in Eq.~\ref{4} has a consistent sequence of roots $\hat{\theta}_n.$

\noindent (b) $\sqrt{n}(\hat{\theta}_n-\theta^g)$ has an asymptotic multivariate normal distribution with mean vector zero and covariance matrix $J_B^{-1}K_B J_B^{-1}$, where $J_B=J_{B}(\theta^g)$, $K_B=K_{B}(\theta^g)$, $K_B(\theta)=Var_g(u_\theta(X)f_\theta(X)B^{\prime\prime}(f_\theta(X)))$.\\

When $g=f_{\theta}$ for some $\theta\in \Theta$ then the above expressions simplify to 
\begin{equation}
\begin{split}
&J_B=\int u_{\theta} u_{\theta}^T f_{\theta}^{2}B^{\prime\prime}(f_\theta),\; K_B=\int u_{\theta} u_{\theta}^T f_{\theta}^{3}(B^{\prime\prime}(f_\theta))^2-\zeta_B\zeta_B^T,\\
& \zeta_B=\int u_{\theta}f_{\theta}^{2}B^{\prime\prime}(f_\theta).
\end{split}
\label{14}
\end{equation}
\label{thm2}
\end{thm}
We are now going to establish that the DPD belongs to the class of $\phi$-DPD. We will also show that under certain conditions a judicial choice of $\phi(\cdot)$ yields estimators which may fit with our aims. Now our unbiased estimating equation for $\phi$-DPD is
\begin{equation}
\begin{split}
&\frac{1}{n}\sum_{i=1}^{n} u_\theta(X_i)f_\theta^{1+\beta}(X_i)\phi(f_\theta(X_i),\gamma)-\\
& \int u_\theta(x) f^{2+\beta}_\theta(x)\phi(f_\theta(x),\gamma)dx=0,
\end{split}
\label{15}
\end{equation}
and under assumptions similar to (A1)-(A5) and $g=f_\theta$, the expressions in Eq.~\ref{14} simplify to
\begin{equation}
\begin{split}
&J_\phi=\int u_{\theta} u_{\theta}^T f_{\theta}^{2+\beta}\phi(f_\theta,\gamma),\\
& K_\phi=\int u_{\theta} u_{\theta}^T f_{\theta}^{3+2\beta}\phi^2(f_\theta,\gamma)-\zeta_\phi\zeta_\phi^T,\\ & \zeta_\phi=\int u_{\theta}f_{\theta}^{2+\beta}\phi(f_\theta,\gamma).
\end{split}
\label{16}
\end{equation} 
A straightforward simplification of the expressions in part (b) of Theorem \ref{thm2} under $\phi$-DPD leads to the general expressions
\begin{equation}
\begin{split}
& K_\phi=\int u_{\theta} u_{\theta}^T f_{\theta}^{2+2\beta}\phi^2(f_\theta,\gamma)g-\zeta_\phi\zeta_\phi^T,\\ 
& \zeta_\phi=\int u_{\theta}f_{\theta}^{1+\beta}\phi(f_\theta,\gamma)g,
\end{split}
\label{17}
\end{equation}
and

\begin{equation}
\begin{split}
J_\phi=&\int u_\theta u_\theta^T f_\theta^{2+\beta}\phi(f_\theta,\gamma)\\
&+\int \Big(\kappa_\theta-\beta u_\theta u_\theta^T\Big)(g-f_\theta)f_\theta^{1+\beta}\phi(f_\theta,\gamma),
\end{split}
\label{18}
\end{equation}

\begin{equation}
\begin{split}
&\kappa_\theta = -\Big[\nabla_\theta u_\theta+\Big(1+\frac{f_\theta \phi^{\prime}(f_\theta,\gamma)}{\phi(f_\theta,\gamma)}\Big)u_\theta u_\theta^T\Big],\\ &\phi^{\prime}(x,\gamma)=\frac{\partial\phi(x,\gamma)}{\partial x}.
\end{split}
\end{equation}

\begin{remark}
Notice that $1+\frac{x\phi^{\prime}(x,\gamma)}{\phi(x,\gamma)}=\frac{1}{\phi(x,\gamma)}\frac{\partial }{\partial x}\left[x\phi(x,\gamma)\right]$. If $\lim_{\gamma\rightarrow 0^+}$ and $\frac{\partial}{\partial x}$ are interchangeable for $\phi(\cdot)$ then by (P4) it can be concluded that $\kappa_\theta$ converges to $i_\theta$ as $\gamma\rightarrow 0^+$.
\end{remark}   

\begin{thm}
If $u_\theta(x) f_\theta(x)^{1+\beta}$, $u_{\theta}(x) u_{\theta}(x)^T f_{\theta}(x)^{1+2\beta}$, $u_{\theta}(x) u_{\theta}(x)^T f_{\theta}(x)^{1+\beta}$ are integrable and $f_\theta(x)\phi(f_\theta(x),\gamma)$ is bounded by some universal constant then the following hold.\\
(a) The usual \textsc{DPD} defined in Eq.~\ref{5} is a special limiting case of $\phi$-\textsc{DPD}.

\noindent (b) If $g=f_\theta$ for some $\theta\in \Theta$ and if for the \textsc{DPD} there exists $ \alpha,\;\beta$ such that the asymptotic relative efficiency (ARE) of the estimator under tuning parameter $\beta$ is greater than that of the estimator under tuning parameter $\alpha$, then there exists $\gamma$ such that $\phi$-\textsc{DPD} with tuning parameter $(\beta,\gamma)$ generates an estimator with higher ARE than the MDPDE with tuning parameter $\alpha$.  
\label{thm3}
 
\end{thm}
If $\phi(x,\gamma)=\frac{1}{\gamma}\log(1+\frac{\gamma}{x})$, then for $x>0$; $x\phi(x,\gamma)<1$ (as $\log(1+y)<y$ for $y>0$) and $u_\theta(x) f_\theta(x)^{1+\beta}$, $u_{\theta}(x) u_{\theta}(x)^T f_{\theta}(x)^{1+2\beta}$, $u_{\theta}(x) u_{\theta}(x)^T f_{\theta}(x)^{1+\beta}$ are integrable under standard parametric models. So Theorem \ref{thm3} holds for such a choice of $\phi(\cdot)$. Symbolically, the divergence generated by the $B(\cdot)$ function obtained through this formulation will be referred to as the logarithmic $\phi$-DPD (or ${\rm L}\phi{\rm DPD}$). We will denote this divergence between the densities $g$ and $f$ corresponding to tuning parameters $\beta$ and $\gamma$ as L$\phi$DPD$_{\beta, \gamma}(g, f)$.

Our choice for $B^{\prime\prime}(\cdot)$ in the ${\rm L}\phi{\rm DPD}$ case is $B^{\prime\prime}(x)=\frac{1}{\gamma}x^\beta\log\left(1+\frac{\gamma}{x}\right) \hspace{0.2cm}  0<\beta\leq 1,\; 0<\gamma\leq 1,\; x>0.$ The corresponding $B$ function may be expressed in the integral form as
\begin{equation}
B(x)=\frac{1}{\gamma}\int_{0}^{x}\int_{0}^{t}s^\beta\log\left(1+\frac{\gamma}{s}\right)\;ds\;dt.
\label{20}
\end{equation}

Obviously other choices are possible, but we have found the L$\phi$DPD to be a very useful divergence for our purpose, and for the rest of the paper all our illustrations will be in relation to the L$\phi$DPD. We will refer to the corresponding minimum distance estimator as ML$\phi$DE.

\subsection{The Influence Function of ML\texorpdfstring{$\phi$}DE}
It is easy to see that the ML$\phi$DE is also an M-estimator. Let the minimum ${\rm L}\phi{\rm DPD}$ functional $T_{\beta,\gamma}(G)$ be defined as 
$$\textsc{${\rm L}\phi{\rm DPD}$}_{\beta,\gamma}(g,f_{T_{\beta,\gamma}(G)})=\inf_{\theta\in \Theta}\textsc{${\rm L}\phi{\rm DPD}$}_{\beta,\gamma}(g,f_\theta).
$$ 
Under $g=f_\theta$ the influence function of this minimum distance estimator simplifies to
\begin{equation*}
\begin{split}
    &\textsc{IF}(y,T_{\beta,\gamma},G)\\
    &=\Big[\int u_\theta u_\theta^T f_\theta^{2+\beta}\log\Big(1+\frac{\gamma}{f_\theta}\Big)\Big]^{-1}\\
    &\Big[u_\theta(y)f_\theta^{1+\beta}(y)\log\Big(1+\frac{\gamma}{f_\theta(y)}\Big)-\int u_\theta f_\theta^{2+\beta}\log\Big(1+\frac{\gamma}{f_\theta}\Big)\Big]
    \end{split}
\end{equation*}
for $0<\beta\leq 1, 0<\gamma\leq 1$. If $\int u_\theta u_\theta^T f_\theta^{2+\beta}\log\Big(1+\frac{\gamma}{f_\theta}\Big)$ and $\int u_\theta f_\theta^{2+\beta}\log\Big(1+\frac{\gamma}{f_\theta}\Big)$ are finite then the expressions in Eq.~\ref{17} are finite if $u_\theta(y)f_\theta^{1+\beta}(y)\log\Big(1+\frac{\gamma}{f_\theta(y)}\Big)$ is finite which is indeed the case for most parametric models suggesting the observed robustness of the ML$\phi$DE under those parametric models.
 \begin{figure}[t]
\centering
\includegraphics[scale=.4]{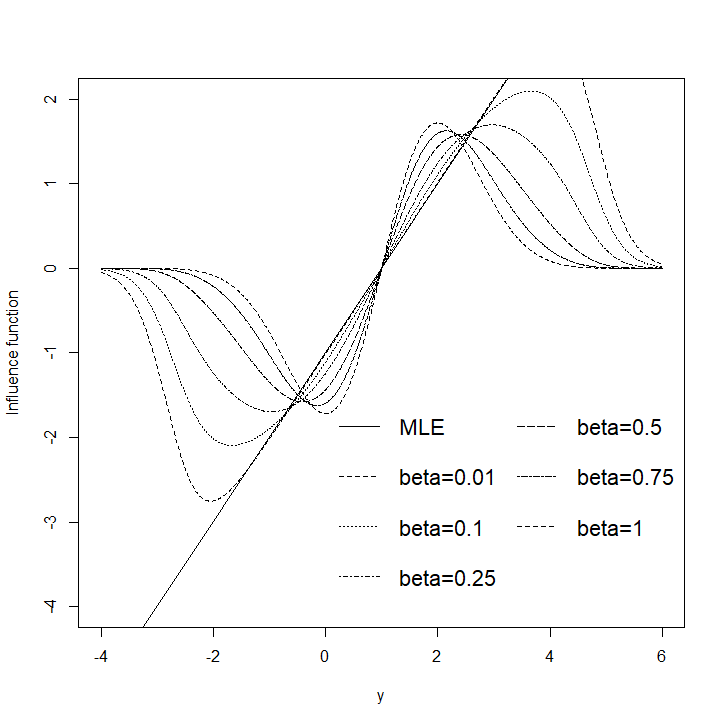}
\caption{Influence function of the ML$\phi$DE for various values of $\beta$ with fixed $\gamma=0.001$ under $N(\mu,1)$ model at $N(1,1)$}
\end{figure}

In Figure 1 it is clearly seen that the tuning parameter $\beta$ has a significant impact on the  robustness of the estimator and the influence functions redescend faster for larger values of $\beta$. On the other hand, for fixed $\beta$ the influence functions are somewhat closer for different $\gamma$ as seen in Figure 2. It suggests that $\gamma$ has a less pronounced impact on robustness than $\beta$, although the graphs in Figure 2 indicate that larger $\gamma$ lead to relatively stronger downweighting.

\begin{figure}[t]
\centering
\includegraphics[scale=.45]{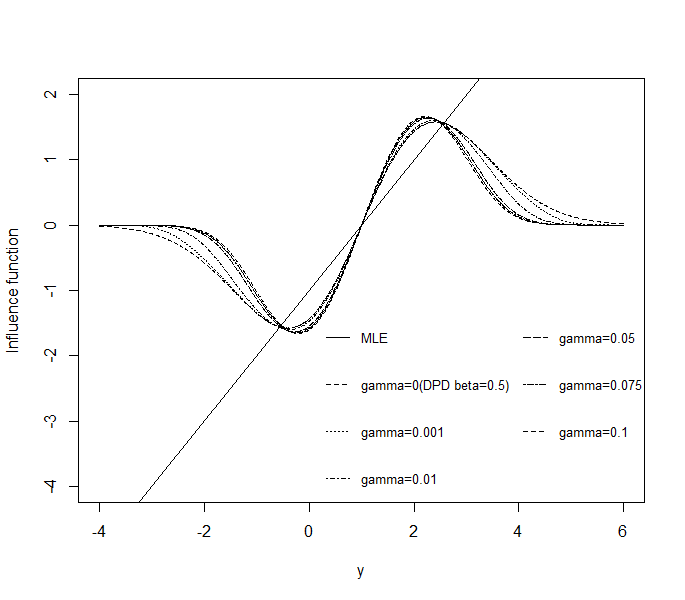}
\caption{Influence function of the ML$\phi$DEs for various values of $\gamma$ with fixed $\beta=0.5$ under the $N(\mu,1)$ model at $N(1,1)$}
\end{figure}
\subsection{The Breakdown Point under the Location Model}
Now we will establish the breakdown point of the minimum ${\rm L}\phi{\rm DPD}$ functional under the location family of densities $\mathcal{F}=\{f_\theta(x)=f(x-\theta):\theta\in\Theta\}$.~Let $B(\cdot)$ be the function defined in Eq.~\ref{20}. Define the quantities:
$$\int B(\epsilon f(x-\theta))dx=\int B(\epsilon f(x))dx:=M_{f,\epsilon}^{(1)},$$
\begin{align*}
    &\int \Big[B(f(x-\theta))+(\epsilon-1)f(x-\theta) B^{\prime}(f(x-\theta))\Big]dx\\
    & =\int \Big[B(f(x))+(\epsilon-1)f(x) B^{\prime}(f(x))\Big]dx \\
    &:=M_{f,(\epsilon-1)}^{(2)}.
\end{align*}
Define $d(g,f)=B(g)-B(f)-(g-f)B^{\prime}(f)$ and let $D(g,f)=\int d(g,f)$.~From Eq.~\ref{20} we have $d(g,0)=\lim_{f\rightarrow 0^+}d(g,f)= B(g)$.

Consider the contamination model $H_{\epsilon,n}=(1-\epsilon)G+\epsilon K_n$, where $\{K_n\}$ is a sequence of contaminating distributions. Let $h_{\epsilon,n}, g$ and $k_n$ be the corresponding densities. We say that there is breakdown in the minimum ${\rm L}\phi{\rm DPD}$ functional for $\epsilon$ level contamination if there exists a sequence $K_n$ such that $|T_{\beta,\gamma}(H_{\epsilon,n})-T_{\beta,\gamma}(G)|\rightarrow \infty$ as $n\rightarrow \infty$. We write below $\theta_n=T_{\beta,\gamma}(H_{\epsilon,n})$ and assume that the true distribution belongs to the model family, i.e., $g=f_{\theta^g}$. We make the following assumptions.

\vspace{0.1in}
\noindent (BP1) $\int \min\{f_\theta(x),k_n(x)\}dx\rightarrow 0$ as $n\rightarrow \infty$ uniformly for $|\theta|\leq c$ for any fixed $c$, i.e., the contamination distribution is asymptotically singular to the true distribution and to specified models within the parametric family.

\vspace{0.1in}
\noindent (BP2) $\int \min\{f_{\theta^g}(x),f_{\theta_n}(x)\}dx\rightarrow 0$ as $n\rightarrow \infty$ if $|\theta_n|\rightarrow \infty$ as $n\rightarrow \infty$, i.e., large values of $\theta$ give distributions which become asymptotically singular to the true distribution.

\vspace{0.1in}
\noindent (BP3) The contaminating sequence $\{k_n\}$ is such that 
$$D(\epsilon k_n,f_\theta)\geq D(\epsilon f_\theta,f_\theta)=M_{f,\epsilon}^{(1)}-M_{f,(\epsilon-1)}^{(2)}$$
for any $\theta\in \Theta$ and $0<\epsilon<1$ and $\limsup_{n\rightarrow\infty}\int B(\epsilon k_n)\leq M_{f,\epsilon}^{(1)}.$
\begin{thm}
Under the assumptions (BP1)-(BP3) above, the asymptotic breakdown point $\epsilon^*$ of the ${\rm L}\phi{\rm DPD}$ functional is at least 0.5 at the location model.
\label{thm4}
\end{thm}
\section{Simulation Study Under \texorpdfstring{$\rm L\phi\rm DPD$} and the Advantages of \texorpdfstring{$\rm ML\phi\rm DE$}{Lg}}
\subsection{Description and Results}
Here we have performed a simulation study to analyze the performance of the ${\rm L}\phi{\rm DPD}$ and the associated minimum distance estimators under the $N(\mu,1)$ model at a given level of contamination. In the following study data are generated from two normal mixtures, $0.9N(0,1)+0.1N(5,1)$ and $0.8N(0,1)+0.2N(5,1)$, where $N(0, 1)$ represents the target distribution and the second component is the contamination. The sample size is 50. The empirical MSE for the location model has been calculated by replicating the process 1000 times, evaluating the estimate for each replication and taking average squared error loss against the target value, i.e., $\mu=0$. In Table 1 the theoretical asymptotic relative efficiency of minimum ${\rm L}\phi{\rm DPD}$ estimator and MDPDE is shown for different values of $(\beta,\gamma)$ while in Table 2 and Table 3 the simulated mean square errors are presented under contaminated normal data under two different contamination levels.

    \begin{table*}[t]
    \begin{center}
\scalebox{0.9}{
    \begin{tabular}{l  l  l  l  l  l  l  l  l   l }
    \hline
    $\beta$ & $\gamma=0$ & $\gamma=0.01$ & $\gamma=0.02$ & $\gamma=0.03$ & $\gamma=0.04$ & $\gamma=0.05$ & $\gamma=0.06$ & $\gamma=0.07$ & $\gamma=0.08$ \\ \hline
    \hline
    0.1 & 98.8 & 95.7 & 94.0 & 92.8 & 91.7 & 90.8 & 90 & 89.4 & 88.7
     \\ 
   0.2& 95.9& 92.6 & 90.8 & 89.5 & 88.4 & 87.4 & 86.6 & 85.9 & 85.3\\ 
   
    0.3 & 92.1 & 89 & 87.2 & 85.9 & 84.8 & 83.9 & 83 & 82.3 & 81.7 \\
    
    0.4 & 88 & 85.2 & 83.5 & 82.2 & 81.1 & 80.2 & 79.4 & 78.7 & 78.1 \\
    
    0.5 & 83.8 & 81.3 & 79.7 & 78.5 & 77.5 & 76.6 & 75.9 & 75.2 & 74.6 \\
    
   0.6 & 79.7 &  77.4 & 76 & 74.9 & 74 & 73.2 & 72.5 & 71.8 & 71.2 \\
   
   0.7 & 75.7 & 73.8 & 72.5 & 71.4 & 70.6 & 69.8 & 69.1 & 68.5 & 68 \\
   
   0.8 & 71.9 & 70.2 & 69 & 68.1 & 67.3 & 66.6 & 66 & 65.4 & 64.9 \\
   
   0.9 & 68.3 & 66.9 & 65.8 & 64.9 & 64.2 & 63.6 & 63 & 62.5 & 62 \\
   
   1 & 65 & 63.7 & 62.7 & 61.9 & 61.2 & 60.7 & 60.1 & 59.7 & 59.2 \\
    \hline

    \end{tabular}}
    \captionof{table}{Asymptotic relative efficiency of the ${\rm ML}\phi{\rm DE}$ and DPDE (\%) for different $(\beta,\gamma)$ under $N(0,1)$ location model. Here the $\gamma=0$ column represents the MDPDE.}
    
    \end{center}
    \end{table*}
    
    \begin{table*}[t]
    \begin{center}
\scalebox{0.9}{
    \begin{tabular}{l  l  l  l  l  l  l  l  l  l }
    
    \hline
    $\beta$ & $\gamma=0$ & $\gamma=0.01$ & $\gamma=0.02$ & $\gamma=0.03$ & $\gamma=0.04$ & $\gamma=0.05$ & $\gamma=0.06$ & $\gamma=0.07$ & $\gamma=0.08$ \\ \hline
    \hline
    $0.1$ & 0.1 & $0.0293$ & $0.0278$ &    $0.029$ &    $0.0259$ &    $0.0271$ &    $0.0278$    & $0.0282$    & $0.0259$
     \\ 
   $0.2$ & 0.056 & $0.0273$ &    $0.0277$ &    $0.0254$    & $0.0252$    & $0.0248$    & $0.0268$    & $0.026$    & $0.0266$
    \\ 
    $0.3$ &  0.036 & $0.0281$ & $0.0257$ &    $0.0267$    & $0.0277$    & $0.0264$    & $0.0279$    & $0.0268$ & $0.0266$
     \\
    $0.4$ & 0.0268 &  $0.0267$ & $0.0261$ & $0.0265$ & $0.027$  & $0.0261$ & $0.0289$ & $0.0273$ & $0.0265$ \\
    
    $0.5$ & 0.0294 & $0.0277$ & $0.0276$ & $0.0276$ & $0.0276$ & $0.0307$ & $0.0284$ & $0.0296$ & $0.0291$ \\
   
       $0.6$ & 0.0275 & $0.0272$ & $0.0298$ & $0.0307$ & $0.0293$ & $0.0305$ & $0.0295$ & $0.0293$ & $0.0296$ \\
    
    $0.7$ & 0.0277 & $0.0276$ & $0.0296$ & $0.0294$ & $0.0301$ & $0.032$ & $0.0311$ & $0.0311$ & $0.03$ \\
    
    $0.8$ & 0.0292 & $0.0305$ & $0.0327$ & $0.0308$ & $0.0342$ & $0.03$ & $0.0288$ & $0.0336$ & $0.0315$ \\
    
    $0.9$ & 0.0309 & $0.0299$ & $0.0313$ & $0.03445$ & $0.0309$ & $0.0295$ & $0.0358$ & $0.0326$ & $0.032$ \\
   
    $1$ & 0.0313 & $0.032$ & $0.035$ & $0.0362$ & $0.0369$ & $0.0361$ & $0.03368$ & $0.034$ & $0.0335$ \\
    \hline
    
    \end{tabular}}
    \captionof{table}{Empirical MSE of the ${\rm ML}\phi{\rm DE}$ and DPDE for different values of $(\beta,\gamma)$ under 10\% contaminated data for location model. Here the $\gamma=0$ column represents the MDPDE.}
    
    \end{center}
    \end{table*}

\begin{table*}[t]
    \begin{center}
\scalebox{0.9}{
    \begin{tabular}{l  l  l  l  l  l  l  l  l  l }
    
    \hline
    $\beta$ & $\gamma=0$ & $\gamma=0.01$ & $\gamma=0.02$ & $\gamma=0.03$ & $\gamma=0.04$ & $\gamma=0.05$ & $\gamma=0.06$ & $\gamma=0.07$ & $\gamma=0.08$ \\ \hline
    \hline
    $0.1$ & 0.3214 & 0.0329 & 0.0314 & 0.0308 & 0.0306 & 0.0305 & 0.0304 & 0.0304 & 0.0305
     \\ 
   $0.2$ & 0.0786 & 0.0312 & 0.0306 & 0.0305 & 0.0305 & 0.0305 & 0.0306 & 0.0307 & 0.0306
    \\ 
    $0.3$ & 0.0414 & 0.0308 & 0.0306 & 0.0307 & 0.0309 & 0.0310 & 0.0311  & 0.0309 & 0.0313
     \\
    $0.4$ & 0.0342 & 0.0311 & 0.0312 & 0.0314 & 0.0316  & 0.0318 & 0.0320 & 0.0322 & 0.0322 \\
    
    0.5 & 0.0327 & 0.0317 & 0.0320 & 0.0323 & 0.0326 & 0.0328 & 0.0330 & 0.0323 & 0.0324 \\
   
    0.6 & 0.0329 & 0.0327 & 0.0331 & 0.0334 & 0.0351 & 0.0372 & 0.0226 & 0.0303 & 0.0310 \\
    
    0.7 & 0.0366 & 0.0413 & 0.0409 & 0.0346 & 0.0366 & 0.0364 & 0.0418 & 0.0421 & 0.0423 \\
    
    0.8 & 0.0382 & 0.0388 & 0.0394 & 0.0399 & 0.0403 & 0.0407 & 0.0410 & 0.0408 & 0.0412 \\
    
    0.9 & 0.0424 & 0.0428 & 0.0432 & 0.0436 & 0.0438 & 0.0442 & 0.0445 & 0.0447 & 0.0309 \\
   
    1 & 0.0437 & 0.0442 & 0.0446 & 0.0390 & 0.0293 & 0.0445 & 0.0267 & 0.0452 & 0.0467 \\
    \hline
    
    \end{tabular}}
    \captionof{table}{Empirical MSE of the ${\rm ML}\phi{\rm DE}$ and DPDE for different values of $(\beta,\gamma)$ under 20\% contaminated data for location model. Here the $\gamma=0$ column represents the MDPDE.}
    
    \end{center}
    \end{table*}

 \subsection{The \texorpdfstring{$\rm L\phi\rm DPD$} versus the DPD}

We briefly note our observations as may be evident from Tables 1 and 2. The asymptotic efficiencies of the minimum divergence estimators decrease with increasing $\beta$ and increasing $\gamma$. Note that given an $\alpha \in (0, 1)$, it may be possible to choose $\beta \in (0,\alpha)$ and $\gamma \in (0, 1)$ so that, in relation to our numerical study, ML$\phi$DE$_{\beta, \gamma}$ beats MDPDE$_\alpha$ both in terms of asymptotic model efficiency and the empirical mean square error under contamination. As an illustration, consider MDPDE$_{0.5}$ in the first contaminated model.~The corresponding MSE and asymptotic relative efficiency are 0.0294 and 83.8\% respectively. Now choose the ${\rm L}\phi{\rm DPD}$ parameter $(\beta,\gamma)=(0.3,0.01)$. In this case, the corresponding MSE and efficiency of the ML$\phi$DE are 0.0281 and 89\% respectively. Similarly ML$\phi$DE$_{0.2, 0.04}$ appears to dominate MDPDE$_{0.4}$ both in terms of asymptotic efficiency and empirical mean square error. In fact, for practically all the MDPDEs that are considered in the Tables 1 and 2 (as also in Tables 1 and 3), there exists a better ML$\phi$DE, both in terms of asymptotic model efficiency and obtained mean square error under contamination. In most of these cases there are several $(\beta, \gamma)$ combinations which provide the domination over a given MDPDE.  Tables 2 and 3 also show that the robust minimum distance estimators hold out well against the outliers at both 10 and 20 percent contamination. Simulation results not presented here indicate that the same holds for higher levels of contamination smaller than 1/2, a consequence of the high breakdown point of the method under location models. 

\section{Algorithm for Finding the Optimal \texorpdfstring{$(\beta,\gamma)$}{Lg}}
 The L$\phi$DPD can generate many different kinds of estimators, starting from the most efficient estimator to highly robust estimators. For example, in the limit $\gamma \rightarrow 0$ and $\beta \rightarrow 0$, one gets the likelihood disparity which is minimized by the classical maximum likelihood estimator. On the other hand, relatively larger values of $\beta$ and $\gamma$ lead to estimators with extremely high outlier stability. In a given situation, therefore, it is imperative that one is able to choose the most suitable tuning parameters for that particular case. Here we consider a data driven algorithm for selecting the ``optimal'' tuning parameters $(\beta, \gamma)$ which would provide best compromise for the given situation. For this purpose we modify an approach of Warwick (2002), pp. 78-82, and minimize an empirical version of the asymptotic summed mean square error. The optimization technique is a two stage process. Suppose that the data are generated by a contaminated version of a model distribution, and let $\theta_0$ be the parameter for the model component. Although the data are generated by a contaminated version, the parameter $\theta_0$ of the model component is our target parameter. The spirit of such a set up is described in Warwick and Jones (2005). Let $\theta_{\beta,\gamma} =T_{\beta,\gamma}(G)$ be the corresponding minimum distance functional and $\hat{\theta}_{\beta,\gamma}$ is the solution of the unbiased equation of ${\rm L}\phi{\rm DPD}$ with tuning parameter $(\beta,\gamma)$ based on the data.~The summed mean square error of the minimum ${\rm L}\phi{\rm DPD}$ estimator has the asymptotic formula
\begin{equation}
\begin{split}
&E\big[\big(\hat{\theta}_{\beta,\gamma}-\theta^*\big)^T\big(\hat{\theta}_{\beta,\gamma}-\theta^*\big)\big]\\
&=\big(\theta_{\beta,\gamma}-\theta
 ^*\big)^T\big(\theta_{\beta,\gamma}-\theta
 ^*\big)+n^{-1} \tr\{var\big(\hat{\theta}_{\beta,\gamma}\big)\}.
 \end{split}
 \label{21}
\end{equation} 
Here $\theta^*$is the pilot estimator playing the role of $\theta_0$ and $\tr\{\cdot\}$ represents the trace of matrix.~The asymptotic covariance matrix of $\sqrt{n}(\hat{\theta}_{\beta,\gamma}-\theta_{\beta,\gamma})$ is $J^{-1}KJ^{-1}$, where $J$ and $K$ are as in Eq.~\ref{16} with $\phi(x,\gamma)=\frac{1}{\gamma}\log(1+\frac{\gamma}{x})$.
So the estimated asymptotic summed mean square of the ML$\phi$DE is
\begin{equation}
\big(\theta_{\beta,\gamma}-\theta
 ^*\big)^T\big(\theta_{\beta,\gamma}-\theta
 ^*\big)+\dfrac{1}{n}J^{-1}KJ^{-1}.
 \label{22}
\end{equation}

For the multiparameter case, the above quantity is a matrix. So trace of the matrix is used to provide a global measure of the summed mean square error for minimization. Thus when there are two parameters to be estimated (say $(\mu,\sigma)$ for $N(\mu,\sigma)$ model) then the expression to be minimized is 
\begin{equation}
\begin{split}
n^{-1}\tr\{
J^{-1}({\theta}_{\beta, \gamma} )
&K({\theta}_{\beta, \gamma} ))
J^{-1}({\theta}_{\beta, \gamma} )\}\\
&+(\mu-\mu^*)^2+(\sigma-\sigma^*)^2.
\end{split}
\label{23}
\end{equation}
The optimal value of $(\beta,\gamma)$ is the minimizer of Eq.~\ref{23} under certain conditions. One important note is that in the first stage of minimization our pilot estimate for $\theta^*$ is taken to be a good robust estimate based on the data as suggested in \cite{Warwick}. The empirical summed mean square error is then obtained by evaluating the expressions in Eq.~\ref{22}  or Eq.~\ref{23} after substituting $\hat{\theta}_{\beta, \gamma}$ for $\theta_{\beta, \gamma}$ and the empirical distribution $G_n$ in place of the true unknown distribution $G$. Let us denote this empirical summed mean square error by AMSE in the following.\\
\\
\underline{\textbf{Algorithm:}}\\

Given a dataset $\textbf{X}_{n\times 1}$ we perform the following steps to obtain the estimate of $\theta$. 
\begin{enumerate}
\item Apply the method suggested in \cite{Warwick} to get an optimal $\alpha$ for MDPDE. Suppose this value is $\alpha_w$. This step is the 1st stage of optimization by assuming an initial pilot estimate of $\theta^*$. 
\item Consider the interval $(0,\alpha_w)$. Update the pilot estimate for $\theta^*= \hat{\theta}_{\alpha_w}$, which is MDPDE of $\theta$ with $\alpha_w$ as the tuning parameter.
\item Perform a two dimensional optimization which selects the value of $(\beta, \gamma)$ for which the minimum
\begin{equation}
\min_{\beta\in(0,\alpha_w)}\Big[\min_{\gamma\in(0,1]}\textsc{AMSE}(\hat{\theta}_{\beta,\gamma})\Big]
\end{equation} 
is attained under the constraint AMSE($\hat{\theta}_{\beta,\gamma})<$ AMSE($\hat{\theta}_{\alpha_w}$).

An alternative to this approach could be to perform an unrestricted minimization of AMSE($\hat{\theta}_{\beta,\gamma})$ with respect to ($\beta,\gamma$) over the set $(0,1)\times(0,1)$.
\end{enumerate}

 \section{Real Data Examples}
 Here we take some real data sets and use our algorithm to find the optimal tuning parameters to be used in estimating the parameters of the model. We worked with two data sets, Newcomb's light speed data and Short's parallax of the sun data, under normality assumptions. We have used the minimum $L_2$ distance estimates as our pilot estimates of $(\mu,\sigma)$.

 \subsection{Newcomb's Data (Speed of Light)}

 This example involves Newcomb’s light speed data (Stigler, 1977, Table 5). The data size is $n=66$. Under the normal model, the MLE of the mean and standard deviation for these data are found to be equal to $26.212$ and $10.664$, respectively. We employ our algorithm for tuning parameter selection and Table 4 reports the optimal tuning  parameters for DPD and L$\phi$DPD, as well as the parameter estimates at these optimal values.~The estimators are extremely close, but the estimated asympmtotic summed mean square, for whatever it is worth, is lower in case of the ML$\phi$DE.\\

 \begin{table}[b]
\begin{center}
\scalebox{0.9}{
\begin{tabular}{c c c}
 \hline
 Category & MDPDE & ${\rm ML}\phi{\rm DE}$\\
 \hline
 Optimal Tuning Parameter & $\alpha=0.3$ & $(\beta,\gamma)=(0.1,0.03)$\\
 Estimate of $\mu$ & 27.62 & 27.57\\
  Estimate of $\sigma$ & 5.01 & 4.93\\
 AMSE & $0.7$ & $0.64$\\
 \hline
 \end{tabular}}
 \caption{Parameter estimates: Newcomb's light speed data.}
\end{center}
\label{tab 1}
\end{table}

\begin{figure}[t]
\centering
\includegraphics[scale=.45]{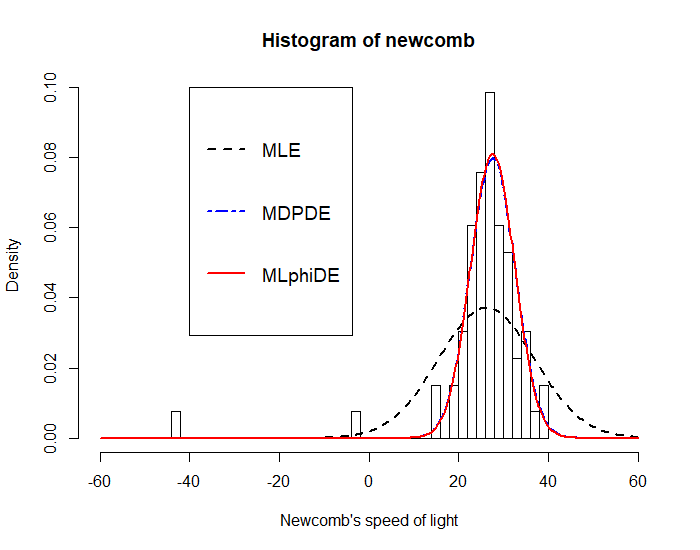}
\caption{Normal density fits for Newcomb's data}
\label{Fig 3}
\end{figure}

\subsection{Short's Data (Parallax of the Sun)}
 This example involves Short’s data for the determination of the parallax of the sun, the angle subtended by the earth’s radius as if viewed and measured from the surface of the sun. From this angle and available knowledge of the physical dimensions of the earth, the mean distance from earth to the sun can be easily determined. The raw observations are presented in Table 4 of Stigler (1977). The data size is $n=53$. Under the normal model, the MLE of the mean and standard deviation for these data are found to be equal to $8.378$ and $0.846$ respectively. We perform all the steps of the aforesaid tuning parameter selection algorithm, and the results of the analysis are now listed in Table 5. Again, the empirical asymptotic MSE for the ML$\phi$DE is slightly better than that of the MDPDE.
 
\vspace{0.3cm}
\begin{table}[b]
\begin{center}
\scalebox{0.9}{
\begin{tabular}{c c c}
 \hline
 Category & MDPDE & ${\rm ML}\phi{\rm DE}$\\
 \hline
 Optimal Tuning Parameter & $\alpha=0.96$ & $(\beta,\gamma)=(0.55,1)$\\
 Estimate of $\mu$ & $8.477$ & $8.478$\\
 Estimate of $\sigma$ & $0.365$ & $0.365$\\
 AMSE & $0.0058$ & 0.0057\\
 \hline
 \end{tabular}}
  \caption{Parameter estimates: Short's data}
\end{center}
\end{table}

\begin{figure}[t]
\centering
\includegraphics[scale=.45]{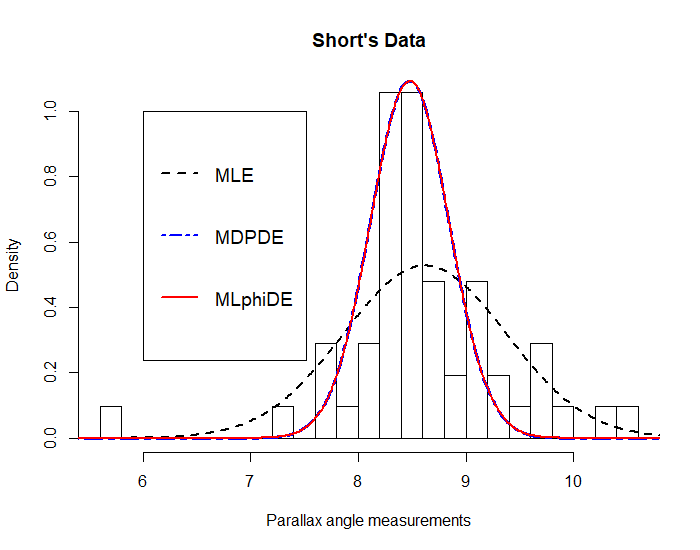}
\caption{Normal density fits for Short's data }
\label{Fig 4}
\end{figure}

From Figure 3 and Figure 4, it is evident that the normal fits coming from the MDPDE and ML$\phi$DE are in the same ballpark. However, if the empirical asymptotic summed mean square error is accepted as a reasonable criterion for discrimination, then the performance of the ML$\phi$DE is better than that of the MDPDE, although the order of improvement is small.

\section{The ML\texorpdfstring{$\phi$}DE for Independent Non-homogeneous Observations}
Here we generalize the above concept to the case of independent but not identically distributed observations. \cite{ghosh2013robust} explains the methodology for this problem in the case of DPD, but here we will extend it to the case of L$\phi$DPD.

Let us assume that the observed data $Y_1, . . . , Y_n$ are independent but for each $i$, $Y_i \sim g_i$ where the densities $g_1, . . . , g_n$ may not be same. We want to model $g_i$ by the family $\mathcal{F}_{i,\theta} = \{f_i(\cdot ; \theta)|\theta \in \Theta\}$ for all $i = 1, 2, . . ., n$. We want to estimate $\theta$ by minimizing the ${\rm L}\phi{\rm DPD}$ between the data and the model. However, the model density may not be same for each $Y_i$'s, and hence we need to calculate the divergence between data and model separately for each data point. For this purpose, we minimize the average divergence between the data points and the models. Therefore, we minimize
$$\dfrac{1}{n}\sum_{i=1}^n d(\hat{g}_i,f_i(\cdot;\theta))$$
with respect to $\theta \in \Theta$, where $d(\hat{g}_i,f_i(\cdot;\theta))$ denotes the
${\rm L}\phi{\rm DPD}$ between the density estimate corresponding to the $i$-th data point and the associated model density. In the presence of only one data point $Y_i$ from density $g_i$, the best possible density estimate of $g_i$ is the (degenerate) density which puts the entire mass on $Y_i$ so that we have
\begin{equation*}
\begin{split}
&d(\hat{g}_i,f_i(\cdot;\theta))\\
&= \frac{1}{\gamma} \int\Big[f_i(y;\theta)\int_0^{f_i(y;\theta)} s^\beta\log\left(1+\frac{\gamma}{s}\right)ds\\ &-\int_0^{f_i(y;\theta)}\int_0^t s^\beta\log\left(1+\frac{\gamma}{s}\right)ds \hspace{0.1cm}dt\Big]dy\\
&-\frac{1}{\gamma} \int_0^{f_i(Y_i;\theta)}s^\beta\log\left(1+\frac{\gamma}{s}\right)ds + K.
\end{split}
\end{equation*}
where $K$ is a constant independent of $\theta$, the parameter of interest. Thus, for the purpose of estimation it suffices to minimize the objective function
\begin{equation}\label{eq:nhmesteq}
H_n(\theta)=\frac{1}{n}\sum_{i=1}^{n} V_i(Y_i;\theta),
\end{equation}
where 
\begin{equation}
\begin{split}
&V_i(Y_i;\theta)\\
&=\frac{1}{\gamma}\int\Big[f_i(y;\theta)\int_0^{f_i(y;\theta)}s^\beta\log\left(1+\frac{\gamma}{s}\right)ds\\
&-\int_0^{f_i(y;\theta)}\int_0^t s^\beta\log\left(1+\frac{\gamma}{s}\right)ds \hspace{0.1cm}dt\Big]dy\hspace{0.1cm}\\
&-\hspace{0.1cm}\frac{1}{\gamma}\int_0^{f_i(Y_i;\theta)}s^\beta\log\left(1+\frac{\gamma}{s}\right)ds.
\end{split}
\end{equation}

Differentiating the above with respect to $\theta$ we get the estimating equation of the minimum ${\rm L}\phi{\rm DPD}$ estimator for non-homogeneous observations as
\begin{equation}
\begin{split}
&\frac{1}{n}\sum_{i=1}^n\Big[f_i(Y_i;\theta)^{\beta+1}\log\left(1+\frac{\gamma}{f_i(Y_i;\theta)}\right)u_i(Y_i;\theta)- \\
&\int f_i(y;\theta)^{\beta+2}\log\left(1+\frac{\gamma}{f_i(y;\theta)}\right)u_i(y;\theta)dy\Big]=0,
\end{split}\label{eq:nhmest}
\end{equation}

where $u_i(\cdot)$ is the score function for $f_i(\cdot)$.

\subsection{Asymptotic Properties}
We will now derive the asymptotic distribution of the minimum ${\rm L}\phi{\rm DPD}$ estimator $\hat{\theta}_n$ defined by the relation
$$H_n(\hat{\theta}_n) = \min_{\theta \in \Theta} H_n(\theta)$$
provided such a minimum exists. Let us first present the necessary set up and conditions. Let the parametric model $\mathcal{F}_{i,\theta}$ be as defined above. We also assume that there exists a best fitting parameter of $\theta$ which is independent of the index $i$ of the different densities. Let us denote it by $\theta^g$. The assumptions hold if all the true densities $g_i$ belong to the model family so that $g_i = f_i(\cdot;\theta)$ for some common $\theta$, and in that case the best fitting parameter is nothing but the true parameter $\theta$.

Next, recall that the ML$\phi$DE $\hat{\theta}_n$ is obtained as a solution of the estimating Eq.~\ref{eq:nhmest}. This equation is satisfied by the minimizer of $H_n(\theta)$ in Eq.~\ref{eq:nhmesteq}. Similarly, we also define, for $i=1,2, \cdots,$
\begin{equation}
\begin{split}
&H^{(i)}(\theta)\\
&=\frac{1}{\gamma}\int\Big[f_i(y;\theta)\int_0^{f_i(y;\theta)}s^\beta\log\left(1+\frac{\gamma}{s}\right)ds\\
&-\int_0^{f_i(y;\theta)}\int_0^ts^\beta\log\left(1+\frac{\gamma}{s}\right)ds \hspace{0.1cm}dt\Big]dy\hspace{0.1cm}\\
&-\hspace{0.1cm}\frac{1}{\gamma}\int\Big[\int_0^{f_i(y;\theta)}s^\beta\log\left(1+\frac{\gamma}{s}\right)ds \Big]g_i(y)dy.
\end{split}
\end{equation}
Note, at the best fitting parameter $\theta^g$, we must have
$$\nabla H^{(i)}(\theta^g)=0, i = 1,2,\cdots$$
We also define, for each $i=1,2,\cdots$ the $p \times p$ matrix $J^{(i)}$ whose $(k,l)$-th entry is given by

\begin{equation}
J^{(i)}_{kl} = E_{g_i}[\nabla_{kl}V_i(Y_i;\theta)],
\end{equation}
where $\nabla_{kl}$ represents the partial derivative with respect to the indicated components of $\theta$. We further define the quantities
\begin{equation}
\Psi_n = \dfrac{1}{n}\sum_{i=1}^n J^{(i)},
\end{equation}

\begin{equation}
\Omega_n = \dfrac{1}{n}\sum_{i=1}^n Var_{g_i}[\nabla V_i(Y_i;\theta)].
\end{equation}
A simple calculation shows that,

\begin{equation}
\begin{split}
&J^{(i)}\\
= &\frac{1}{\gamma}\int u_i(y;\theta^g) u_i^T(y;\theta^g)f_i^{\beta+2}(y;\theta^g) \log\Big(1+\frac{\gamma}{f_i(y;\theta^g)}\Big)dy\\
&- \frac{1}{\gamma}\int \Big[\{\nabla u_i(y;\theta^g) + (\beta+1) u_i(y;\theta^g) u_i^T(y;\theta^g)\}\\ &\log\Big(1+\frac{\gamma}{f_i(y;\theta^g)}\Big)
- u_i(y;\theta^g) u_i^T(y;\theta^g)\Big(\frac{\gamma}{\gamma + f_i(y;\theta^g)}\Big) \Big]\\ 
& \{g_i(y;\theta^g) -f_i(y;\theta^g)\} f_i^{\beta+1}(y;\theta^g) dy
\end{split}
\end{equation}
and

\begin{equation}
\begin{split}
&\Omega_n \\
&= \frac{1}{n}\sum_{i=1}^n \frac{1}{\gamma}\int \Big\{\int_0^{f_i(y;\theta)} s^\beta  \log\Big(1+\frac{\gamma}{s}\Big)ds\Big\}^2 g_i(y;\theta) dy\\
&-\frac{1}{n}\sum_{i=1}^{n}\xi_i \xi_i^T ,
\end{split}
\end{equation}
where

\begin{equation}
\xi_i = \frac{1}{\gamma}\int \Big\{\int_0^{f_i(y;\theta)} s^\beta  \log\Big(1+\frac{\gamma}{s}\Big)ds\Big\} g_i(y;\theta) dy.
\end{equation}

We will make the following assumptions to establish the asymptotic properties of the ML$\phi$DE:

\vspace{0.1in}
\noindent (G1) The support $\mathcal{X} = \{y|f_i(y;\theta)>0\}$ is independent of $i$ and $\theta$ for all $i$; the true distributions $G_i$ are also supported on $\mathcal{X}$ for all $i$.

\vspace{0.1in}
\noindent (G2) There is an open subset $\omega$ of the parameter space $\Theta$, containing the best fitting parameter $\theta^g$ such that for almost all $y \in \mathcal{X}$, and all $\theta \in \Theta$, all $i = 1, 2,\cdots$, the density $f_i(y; \theta)$ is thrice differentiable with respect to $\theta$ and the third partial derivatives are continuous with respect to $\theta$.

\vspace{0.1in}
\noindent (G3) For each $i = 1, 2, \cdots$, the three integrals $\int f_i(y;\theta)\int_0^{f_i(y;\theta)} s^\beta\log\left(1+\frac{\gamma}{s}\right)ds\hspace{0.1cm}dy$,\\ $\int\int_0^{f_i(y;\theta)}\int_0^t s^\beta\log\left(1+\frac{\gamma}{s}\right)ds \hspace{0.1cm}dt\hspace{0.1cm}dy$,\\
and $\int\Big[\int_0^{f_i(y;\theta)} s^\beta\log\left(1+\frac{\gamma}{s}\right)ds \Big]g_i(y)dy$
can be differentiated thrice with respect to $\theta$, and the derivatives can be taken under the integral sign (the first indefinite integral).

\vspace{0.1in}
\noindent (G4) For each $i = 1, 2, \cdots$, the matrices $J^{(i)}$ are positive definite and
$$\lambda_0 = \inf_n [\text{min eigenvalue of } \Psi_n] > 0.$$

\vspace{0.1in}
\noindent (G5) There exists functions $M_{jkl}^{(i)}(Y)$ such that
$$|\nabla_{jkl} V_i(Y;\theta)| \leq M_{jkl}^{(i)}(Y) \hspace{0.3cm} \forall \theta \in \Theta, \hspace{0.3cm} \forall i$$

with $E_{g_i}|M_{jkl}^{(i)}(Y)| < \infty \hspace{0.3cm} \forall j,\;k,\;l.$

\vspace{0.1in}
\noindent (G6) For all $j,k$, we have

\begin{equation}\label{33}
\lim_{N \rightarrow \infty} \sup_{n>1} \Big\{\frac{1}{n} \sum_{i=1}^n E_{g_i}[|\nabla_j V_i(Y;\theta)| I(|\nabla_j V_i(Y;\theta)|>N)] \Big\} = 0,
\end{equation}

\begin{equation}\label{34}
\begin{split}
&\lim_{N \rightarrow \infty} \sup_{n>1} \Big\{\frac{1}{n} \sum_{i=1}^n E_{g_i}[|\nabla_{jk} V_i(Y;\theta)- E_{g_i}(\nabla_{jk} V_i(Y;\theta))|\\
& \times I(|\nabla_{jk} V_i(Y;\theta)- E_{g_i}(\nabla_{jk} V_i(Y;\theta))|>N)]\Big\} = 0.
\end{split}
\end{equation}
Here $I(\cdot)$ stands for indicator function.

\vspace{0.1in}
\noindent (G7) For all $\epsilon > 0$, we have

\begin{equation}
\begin{split}
\lim_{n \rightarrow \infty} \Big\{\frac{1}{n} \sum_{i=1}^n E_{g_i}\big[\|\Omega_n^{-1/2}\nabla V_i(Y;\theta)\|^2 & I(\|\Omega_n^{-1/2}\nabla V_i(Y;\theta)\|\\ &> \epsilon \sqrt{n}) \big] \Big\}= 0
\end{split}
\end{equation}

\begin{thm}\label{newthm5}
Under assumptions (G1)-(G7), the following results hold:
\begin{enumerate}[label = (\roman*)]
    \item There exists a consistent sequence $\theta_n$ of roots to the minimum ${\rm L}\phi{\rm DPD}$ estimating Eq.~\ref{eq:nhmest}.
    \item The asymptotic distribution of $\Omega_n^{-\frac{1}{2}} \Psi_n[\sqrt{n}(\theta_n - \theta^g)]$ is $p$-dimensional normal with (vector) mean $0$ and covariance matrix $I_p$, the $p$-dimensional identity matrix.
\end{enumerate}
\end{thm}

Note that, putting $f_i = f$ for all $i$, we get back the corresponding asymptotic properties of the minimum ${\rm L}\phi{\rm DPD}$ estimator for the i.i.d. case. If $f_i = f, i =1, 2,\cdots$, we get $J^{(i)} = J$ for all $i$; thus $\Psi_n=J$ and $\Omega_n = K$. Here $J$ and $K$ are as defined previously. In this case assumptions (G1)–(G5) are exactly the same as the assumptions (A1)-(A5), while assumptions (G6) and (G7) are automatically satisfied by the dominated convergence theorem. Thus the result, which establishes the consistency and asymptotic normality of the minimum ${\rm L}\phi{\rm DPD}$ estimator $\hat{\theta}$ with $n^{1/2}(\hat{\theta} - \theta^g)$ having the asymptotic covariance matrix $\Psi_n^{-1}\Omega_n \Psi_n^{-1} = J^{-1}KJ^{-1}$, emerges as a special case of Theorem \ref{newthm5}.

\subsection{Normal Linear Regression}
A natural situation where the theory proposed above would be immediately applicable is the case of linear regression. We consider the linear regression model
\begin{equation}
y_i = x_i^T \beta + \epsilon_i, \hspace{0.5cm} i = 1,\cdots,n,
\end{equation}
where the error $\epsilon_i$'s are i.i.d. normal variables with mean zero and variance $\sigma^2$, $x_i^T = (x_{i1},\cdots,x_{ip})$ is the vector of the independent variables corresponding to the $i$-th observation and $\beta = (\beta_1,\cdots,\beta_p)^T$ represents the regression coefficients. We will assume that $x_i$'s are fixed. Then $y_i \sim N(x_i^T \beta, \sigma^2)$, and hence the $y_i$’s are independent but not identically distributed. Thus $y_i$’s satisfy our independent but non-homogenous set-up and hence the ML$\phi$DE of the parameter $\theta = (\beta^T, \sigma^2)^T$ can be obtained by minimizing the expression in Eq.~\ref{eq:nhmesteq} with $f_i \equiv N(x_i^T \beta, \sigma^2)$.

\section{Real Data Examples in Regression}
We now consider some real data examples to illustrate the above technique in linear regression. 

\subsection{Hertzsprung-Russel Data}
This example involves a robust regression on the Hertzsprung-Russel data. These data, associated with the Hertzsprung-Russel diagram of the star cluster CYG OB1 containing 47 stars in the direction of Cygnus has been analyzed previously by several authors including \cite{rousseeuw1987robust}.\\

We fit the simple linear regression model $y=\eta_0+\eta_1x+\epsilon$ under homoscedastic normal errors. Here the independent variable $(x)$ is the logarithm of the temperature of the stars, and the dependent variable $(y)$ is the logarithm of the light intensity of the stars. The initial regression parameter values are the least median of squares (LMS) estimates. The initial scale estimate is the scaled median absolute deviation (MAD) of the LMS residuals. We perform the previously mentioned steps of optimal tuning parameter selection and obtain the estimates for the regression coefficients, which are given in Table 6. The regression lines for LS regression, LMS regression and minimum ${\rm L}\phi{\rm DPD}$ regression are given in the Figure 5. The robust performance of the ML$\phi$DE is self evident.

\vspace{0.3cm}
\begin{table}[b]
\begin{center}
\begin{tabular}{c c}
 \hline
 Category & ${\rm ML}\phi{\rm DE}$\\
 \hline
 Tuning Parameter & $(\beta,\gamma)=(1,0.9)$\\
 Estimate of $\eta_0$ & $-8.5557324$\\
  Estimate of $\eta_1$ & $3.0590795$\\
  Estimate of $\sigma$ & $0.4266284$\\
 \hline
 \end{tabular}
 \caption{Regression estimates for Hertzsprung-Russel data}
\end{center}
\label{tab 3}
\end{table}

 \begin{figure}[t]
\centering
\includegraphics[scale=.47]{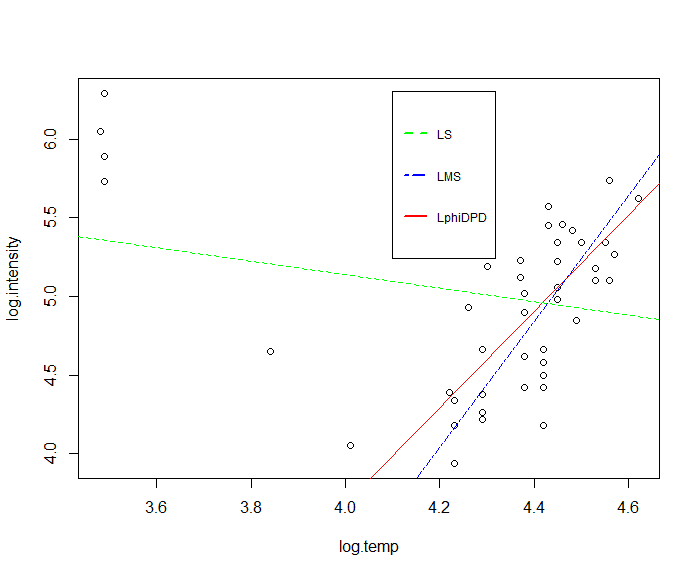}
\caption{Regression fits for the Hertzsprung-Russel data}
\label{Fig 5}
\end{figure}
\subsection{Salinity Data}
This example involves the Salinity data (Table 5, Chapter 3, Rousseeuw and Leroy, 1987). These data were originally presented by \cite{ruppert1980trimmed}. The measurements of the salt concentration of the water and the river discharge taken in North Carolina’s Pamlico Sound were recorded as the data. These data represent a multiple linear model with salinity as the dependent variable $(y)$, and salinity lagged by two weeks $(x_1)$, the number of biweekly periods elapsed since the beginning of the spring season $(x_2)$, and the volume of river discharge into the sound $(x_3)$ as the dependent variable.

We fit the multiple linear regression model $y=\eta_0+\eta_1x_1+\eta_2x_2+\eta_3x_3+\epsilon$ under homoscedastic normal errors. The initial regression parameter values are the least median of squares (LMS) estimates. The initial scale estimate is the scaled median absolute deviation (MAD) of the LMS residuals.

The optimal parameters obtained through our algorithm for optimal parameter selection are presented in Table 7. The residual plots for LS regression, LMS regression and minimum ${\rm L}\phi{\rm DPD}$ regression are given in the Figure 6. Like the LMS method (and unlike the LS method) the ML$\phi$DE gives a nice outlier resistant fit.
\vspace{-0.3cm}
\begin{table}[b]
\begin{center}
\scalebox{0.9}{\begin{tabular}{c c}
 \hline
 Category & ${\rm ML}\phi{\rm DE}$\\
 \hline
 Tuning Parameter & $(\beta,\gamma)=(1,0.9)$\\
 Estimate of $\eta_0$ & $57.16780461$\\
  Estimate of $\eta_1$ & $0.06010002$\\
  Estimate of $\eta_2$ & $-0.01301208$\\
  Estimate of $\eta_3$ & $-2.08372562$\\
  Estimate of $\sigma$ & $0.56157558$\\
 \hline
 \end{tabular}}
 \caption{Regression estimates for Salinity data}
\end{center}
\label{tab 4}
\end{table}

 \begin{figure}[t]
\centering
\includegraphics[scale=.4]{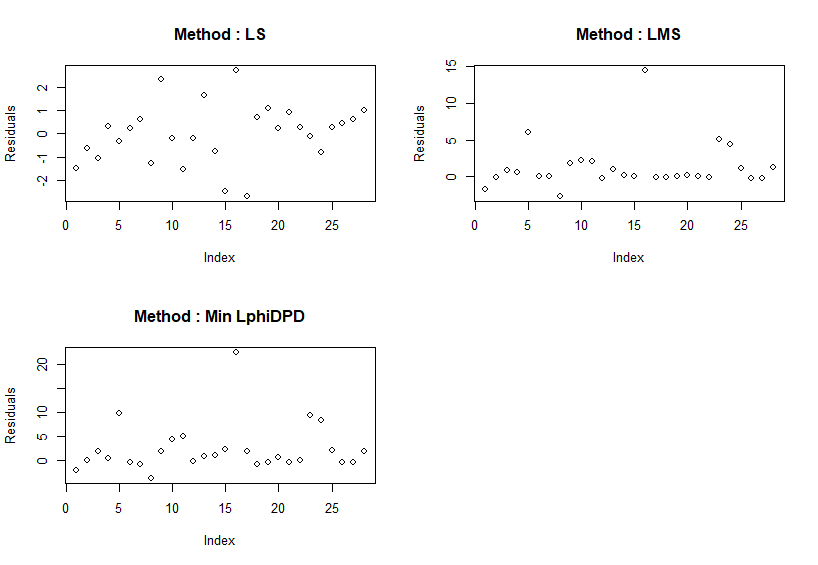}
\caption{Residual plots of the fitted regression models for Salinity data using LS, LMS and minimum ${\rm L}\phi{\rm DPD}$ estimation }
\label{Fig 6}
\end{figure}

\section{Hypothesis Testing using \texorpdfstring{${\rm L}\phi{\rm DPD}$}{Lg}}
Now we develop the tests of parametric hypothesis based on ${\rm L}\phi{\rm DPD}$ divergence. The most common problem is that of testing a simple null hypothesis for a parametric family of densities $\{f_\theta:\theta \in \Theta \subset\mathbb{R}^p\}$ under the one sample case. Here we test
\begin{equation}
    H_0: \hspace{0.1cm}\theta = \theta_0 \quad\text{versus}\quad H_1: \hspace{.1cm} \theta \neq \theta_0
    \label{nullhypo}
\end{equation}
when a random sample $X_1,X_2,\ldots,X_n$ is available from the population of interest. We propose our test statistic as $$T=T_{\beta,\gamma}(\hat{\theta},\theta_0)=2n d_{\beta,\gamma}(f_{\hat{\theta}},f_{\theta_0})$$ where 
\begin{equation}
\begin{split}
&d_{\beta,\gamma}(f_{\hat{\theta}},f_{\theta_0})\\
&=\int \Big[B(f_{\hat{\theta}}(x))-B(f_{\theta_0}(x))-(f_{\hat{\theta}}(x)-f_{\theta_0}(x))B^{\prime}(f_{\theta_0}(x))\Big]dx,
\end{split}
\end{equation}
with $\hat{\theta}=\hat{\theta}_{\beta,\gamma}$ being the ${\rm ML}\phi{\rm DE}$ estimate of $\theta$ and $B(\cdot)$ is as defined in Eq.~\ref{20}.~We shall find the asymptotic distribution of $T$ under $H_0$ and reject the null hypothesis for large values of $T$.

We assume the following regularity conditions of the parametric family of distributions,

\vspace{0.1in}
\noindent (B1) The support of the distribution function $F_{\theta}$, i.e.~ the set $\mathcal{X}=\{x|f_\theta (x)>0\}$ is independent of $\theta$.

\vspace{0.1in}
\noindent (B2) There is an open subset $\omega$ of the parameter space $\Theta$, containing the true parameter value $\theta_0$ such that for almost all $x\in \mathcal{X}$, and all $\theta\in \omega$, the density $f_\theta(x)$ is three times differentiable with respect to $\theta$ and the third partial derivatives are continuous with respect to $\theta$. 

\vspace{0.1in}
\noindent (B3) The integrals $\int B^{\prime\prime}(f_\theta(x))f_\theta^2(x)dx$ can be differentiated with respect to $\theta$, and the derivatives can be taken under the integral sign.

\vspace{0.1in}
\noindent (B4) The $p\times p$ matrix $J(\theta)$ defined by
\begin{equation*}
\begin{split}
    &J_{B,kl}(\theta)=\\
    &E_{\theta}\Big\{\nabla_{kl}\Big(\int [B^{\prime}(f_\theta(x))f_{\theta}(x)-B(f_\theta(x))]dx-B^{\prime}(f_\theta(X))\Big)\Big\}
    \end{split}
\end{equation*}
is positive definite where $E_{\theta}$ represents the expectation under the density $f_{\theta}$. 

\vspace{0.1in}
\noindent (B5) There exists functions $M_{jkl}(x)$ with finite expectation, $j,k,l=1,\ldots, p$, such that
\begin{equation*}
   \begin{split} 
 &\Big|\nabla_{jkl}\Big(\int [B^{\prime}(f_\theta(x))f_{\theta}(x)-B(f_\theta(x))]dx-B^{\prime}(f_\theta(X))\Big)\Big|\\
 &\leq M_{jkl}(X);\;\; \forall \theta \in \omega.
 \end{split}
 \end{equation*}

Then we have the following theorem.
\begin{thm}
Under the assumptions (B1)-(B5) and under the null hypothesis $H_0:\hspace{0.1cm}\theta=\theta_0$ the asymptotic distribution of $T_{\beta,\gamma}(\hat{\theta},\theta_0)$ coincides with the distribution of $$\sum_{i=1}^r \lambda_i Z_i^2,$$ where $Z_i$s are independent standard normals and $\lambda_i$'s are non-zero eigenvalues of $A(\theta_0)\Sigma(\theta_0)$ and $$r=rank(\Sigma(\theta_0)A(\theta_0)\Sigma(\theta_0))$$ $$A(\theta_0)=\nabla_\theta^2 d_{\beta,\gamma}(f_{\theta},f_{\theta_0})|_{\theta=\theta_0}$$
where $\Sigma(\theta_0)$ is the asymptotic covariance matrix of $\sqrt{n} \hat{\theta}_{\beta, \gamma}$ under the null hypothesis and $\nabla^2_\theta$ represents second derivative with respect to $\theta$. . 
\label{thm5}
\end{thm}

We can extend this theorem and hence the testing result to the general two sample problem of testing $H_0:\hspace{0.1cm}\theta_1=\theta_2$ against $H_1:\hspace{0.1cm}\theta_1\neq\theta_2$
where there is a random sample of size $n$ from population $1$ with parameter $\theta_1$ and that of size $m$ from population $2$ with parameter $\theta_2$. Let $\hat{\theta}_1$ and $\hat{\theta}_2$ be ML$\phi$DEs of the parameter in populations 1 and 2, respectively. Then under the (B1)-(B5) regularity conditions on the model, we have the following results.
\begin{thm}
Under the null $H_0: \hspace{0.1cm}\theta_1=\theta_2$, the asymptotic distribution of $$S=S_{\beta,\gamma}(\hat{\theta}_1,\hat{\theta}_2)=\frac{2mn}{m+n}d_{\beta,\gamma}(f_{\hat{\theta}_1},f_{\hat{\theta}_2})$$ coincides with that of $$\sum_{i=1}^r \lambda_i Z_i^2$$ where $Z_i$s are independent standard normals and $\lambda_i$'s are non-zero eigenvalues of $A(\theta_1)\Sigma(\theta_1)$ and $r=rank(\Sigma(\theta_1)A(\theta_1)\Sigma(\theta_1))$ where $A(\theta)$ and $\Sigma(\theta)$ are defined in the statement of Theorem \ref{thm5}. 
\label{thm6}
\end{thm}

\subsection{Equivalence with the Score Test}
A score test, developed in the same spirit under the same set up as in Theorem \ref{thm5}, also has the same asymptotic null distribution. 

\begin{thm}
The score test statistic using the ${\rm L}\phi{\rm DPD}$ for testing the simple null in Eq. \ref{nullhypo} can be given by $$n\bar{U}^T(\theta_0)J_B^{-1}(\theta_0)A(\theta_0)J_B^{-1}(\theta_0)\bar{U}(\theta_0)$$ where $$U_\theta(x)=u_{\theta}(x)B^{\prime\prime}(f_{\theta}(x))f_{\theta}(x)-\int u_{\theta}(x)B^{\prime\prime}(f_{\theta}(x))f_{\theta}^2(x)dx$$ and $$\bar{U}(\theta)=\frac{1}{n}\sum_{i=1}^nU_\theta(X_i)$$ with $$B^{\prime\prime}(x)=\frac{x^\beta}{\gamma}\log\left(1+\frac{\gamma}{x}\right)$$ $$J_B(\theta_0)=-E_{\theta_0}\frac{\partial}{\partial \theta}U_\theta(X_1)|_{\theta=\theta_0}$$ and $A(\theta_0)$ is as described in Theorem \ref{thm5}. Under the null hypothesis, the asymptotic distribution of this statistic is same as that of $T_{\beta,\gamma}(\hat{\theta},\theta_0)$.
\label{thm7}
\end{thm}
\subsection{Divergence Difference test statistic}
We assume that we have a parametric model $\mathcal{F}$ of densities and $X_1,\ldots,X_n$ be i.i.d. from the true distribution $G$ with the same support as the distributions in $\mathcal{F}$. Consider the null hypothesis
\begin{equation}
    H_0: \theta\in \Theta_0\quad \text{versus}\quad H_1:\theta\in \Theta\setminus \Theta_0,
    \label{eq 35}
\end{equation}
where $\Theta_0$ is a proper subset of $\Theta$. The likelihood ratio test (LRT) is one of the most common tests that may be employed in this situation. Define $$\lambda=\frac{\sup_{\theta\in \Theta_0}\mathcal{L}(\theta|X_1,\ldots,X_n)}{\sup_{\theta\in \Theta}\mathcal{L}(\theta|X_1,\ldots,X_n)},$$
where $\mathcal{L}(\theta|X_1,\ldots,X_n)$ is the likelihood of $\theta$ given the data. The test statistic in this case is $-2\log\lambda$. Assume that the distribution function $G$ is discrete. In particular let its support be $\mathcal{X}=\{0,1,2,\ldots\}$, which is also the common support of the family $\mathcal{F}$. Then the test statistic can be expressed in terms of observed relative frequencies $\nu_n$ as
\begin{equation}
\begin{split}
&-2\log \lambda\\&=2\left[\log\left(\prod_{i=1}^{n}f_{\hat{\theta}}(X_i)\right)-\log\left(\prod_{i=1}^{n}f_{\hat{\theta}_0}(X_i)\right)\right]\\
    &= 2n \left[\text{LD}(\nu_n,f_{\hat{\theta}_0})-\text{LD}(\nu_n,f_{\hat{\theta}})\right],
    \label{eq 37}
\end{split}
\end{equation}
where $\text{LD}(\cdot\;,\cdot)$ stands for the likelihood disparity. Here $\hat{\theta}$ and $\hat{\theta}_0$ stands for unrestricted maximum likelihood estimator and maximum likelihood estimator under null hypothesis respectively. Eq.~\ref{eq 37} gives a motivation to construct a new test statistic based on L$\phi$DPD. 

As an analog of the likelihood ratio test, we consider the divergence difference test (DDT) based on L$\phi$DPD to test the hypothesis given in Eq.~\ref{eq 35}. Note that the test statistic in Eq.~\ref{eq 37} can be viewed as a difference of the minimized value of likelihood disparity under null and unrestricted minimum of likelihood disparity.~In the same spirit one may define the following test statistic
\begin{equation}
    \text{DDT}_{\beta,\gamma}(\nu_n)=2n\left[d_{\beta,\gamma}(\nu_n,f_{\hat{\theta}_0})-d_{\beta,\gamma}(\nu_n,f_{\hat{\theta}})\right],
    \label{eq 38}
\end{equation}
$\hat{\theta}_0$ and $\hat{\theta}$ are ML$\phi$DE under null hypothesis and unrestricted minimum ML$\phi$DE respectively. Also note that
\begin{equation*}
\begin{split}
&d_{\beta,\gamma}(\nu_n,f_\theta)\\
&=\sum_{x\in \mathcal{X}}\left[B(\nu_n(x))-B(f_\theta(x))-(\nu_n(x)-f_\theta(x))B^\prime(f_\theta(x))\right],
\end{split}
\end{equation*}
where $B(\cdot)$ is defined as Eq. $\ref{20}$. We will show that under certain regularity conditions the asymptotic distribution of the the test statistic $\text{DDT}_{\beta,\gamma}(\nu_n)$ coincides with the distribution of linear combination of independent chi-squared random variables. Suppose that $\Theta_0$ is defined by a set of $r\leq p$ restrictions on $\Theta$ defined by $R_i(\theta)=0,\, 1\leq i\leq r$.~We assume that the parameter space under $H_0$ can be described through a parameter $\xi=(\xi_1,\ldots,\xi_{p-r})$, with $p-r$ independent components, i.e., $H_0$ specifies that there exists a function $b:\mathbb{R}^{p-r}\to\mathbb{R}^p$ where $\theta=b(\gamma)$, $\gamma\in \Gamma\subseteq \mathbb{R}^{p-r}$. The function $b$ is assumed to have continuous derivative $\dot{b}(\xi)$ of order $p\times(p-r)$ with rank $p-r$. Then the constrained estimator is $\hat{\theta}_0=b(\hat{\xi})$, where $\hat{\xi}$ is the  ML$\phi$DE under the $\xi$ formulation of the model. Let $G=F_\theta$ be the true distribution which belongs to the family $\mathcal{F}$ with parameter $\theta$. Under $H_0$, let $\xi$ be the true value of the reduced parameter. So we have $\theta=b(\xi)$. When the null hypothesis is true under standard regularity conditions it can be easily shown that $\hat{\xi}$ and $\hat{\theta}_0$ are consistent for $\xi$ and $\theta$ respectively in the sense that
\begin{equation}
\begin{split}
    \hat{\xi}&= \xi+n^{-1/2}\left[\dot{b}(\xi)^TJ_B(b(\xi))\dot{b}(\xi)\right]^{-1}\dot{b}(\xi)^T Z_n(b(\xi))\\
    &+o_P(n^{-1/2}),
    \end{split}
    \label{eq 39}
\end{equation}
where $Z_n(b(\xi))$ is $AN(0,K_B(b(\xi)))$. Here $J_B(\cdot)$ and $K_B(\cdot)$ is defined as in Theorem \ref{thm2}. Now we will lay out some appropriate regularity conditions under which we will derive the asymptotic distribution of $\text{DDT}_{\beta,\gamma}(\nu_n)$ under the null hypothesis.

\vspace{0.1in}
\noindent(C1) The assumptions (A1)-(A5) hold under the model conditions.

\vspace{0.1in}
\noindent(C2) The unconstrained minimum L$\phi$DPD estimator $\hat{\theta}$ satisfies
\begin{equation}
    \hat{\theta}=\theta+n^{-1/2}J_B^{-1}(\theta)Z_n(\theta)+o_p(n^{-1/2}),
    \label{eq 40}
\end{equation}
where $Z_n(\theta)$ is $AN(0,K_B(\theta))$.

\vspace{0.1in}
\noindent(C3) The null hypothesis $H_0$ is either simple and $\Theta_0=\{\theta_0\}$, where $\theta_0$ is in the interior of $\Theta$, or $H_0$ is composite and $\Theta_0=\{b(\xi): \xi\in \Gamma\subseteq \mathbb{R}^{p-r}\}$.

\vspace{0.1in}
\noindent(C4) If $H_0$ is composite then the constrained estimator $\hat{\theta}_0=b(\hat{\xi})$ and $\hat{\xi}$ satisfies Eq.~\ref{eq 39}.
Define $$\Sigma_{B,b}(\theta,\xi)=\Tilde{J}^{-1}_{B,b}(\theta,\xi)K_B(\theta)\Tilde{J}^{-1}_{B,b}(\theta,\xi),$$
where $$\Tilde{J}_{B,b}(\theta,\xi)=\left[J_B(\theta)^{-1}-\dot{b}(\xi)[\dot{b}(\xi)^T J_B(\theta)\dot{b}(\xi)]^{-1}\dot{b}(\xi)^T\right]^{-1}$$. 
\begin{thm}
Suppose that assumption (C1)-(C4) hold. Under $f_{\theta_0}$, $\theta_0\in \Theta_0$, the limiting distribution of the distance difference test statistic in Eq.~\ref{eq 38} coincides with the distribution of 
$$\sum_{i=1}^{m}\lambda_i Z_i^2,$$
where $\lambda_i$'s are non-zero eigenvalues of $A(\theta_0)\Sigma_{B,b}(\theta_0,\xi)$ and $m=\text{rank}(A(\theta_0)\Sigma_{B,b}(\theta_0,\xi))$. Moreover if $\Theta_0=\{\theta_0\}$ then under the null hypothesis the asymptotic distribution of distance difference test statistic in Eq.~\ref{eq 38} is same as that of    $T_{\beta,\gamma}(\hat{\theta},\theta_0)$ in Theorem \ref{thm5}.
\label{thm 8}
\end{thm}
\begin{remark}
In the above theorems the null distribution of the test statistic turns out to be same as that of  a linear combination of independent chi squared random variables. In general it is hard to get hold of critical values under this distribution for actually performing the test. Also calculations regarding this distribution become numerically hard. This gives the motivation to explore another test statistic which will lead to a simpler null distribution.
\end{remark}
\subsection{Wald Type Test}
Assume a similar setup of hypothesis testing as in Eq.~\ref{eq 35}. Suppose that the null space $\Theta_0\subseteq \Theta\subseteq \mathbb{R}^p$ is defined by a set of $r\leq p$ restrictions on $\Theta$ defined by $R_i(\theta)=0,\, 1\leq i\leq r$.~Let $G=F_\theta$ be the true distribution which belongs to the family $\mathcal{F}$ with parameter $\theta$. Assume $\hat{\theta}$ to be the ML$\phi$DE of the true parameter $\theta$. Define $\boldsymbol{R}(\theta)=(R_1(\theta),\ldots,R_r(\theta))^T$ and $\boldsymbol{D}(\theta)=\left[\frac{\partial R_i(\theta)}{\partial \theta_j}\right]_{r\times p}$.~Under the spirit of the original Wald test statistic, we can construct the following test statistic
\begin{equation*}
    W(\hat{\theta})=\boldsymbol{R}(\hat{\theta})^T\left(\boldsymbol{D}(\hat{\theta})\boldsymbol{\Sigma}(\hat{\theta})\boldsymbol{D}(\hat{\theta})^T\right)^{-1}\boldsymbol{R}(\hat{\theta}),
\end{equation*}
where $\boldsymbol{\Sigma}(\theta)=J_B(\theta)^{-1}K_B(\theta)J_B(\theta)^{-1}$ under the $B(\cdot)$ function described in Eq.~\ref{20}.~Under standard regularity conditions it is easy to prove that the asymptotic distribution of $W(\hat{\theta})$ is $\chi^2_r$ under the null hypothesis. The proof follows from simple application of delta method theorem on the quantity $\boldsymbol{R}(\hat{\theta})$ and the fact that under the null hypothesis $\sqrt{n}(\hat{\theta}-\theta)$ is $AN(0,\boldsymbol{\Sigma}(\theta))$. The main benefit of this test statistic is that its asymptotic null distribution is simpler.~Hence it is easy to perform numerical computations based on these statistics. For example, the critical values of the test statistic can be computed with ease in this case.
\subsection{Real Data Example}
Researchers needed to evaluate the effectiveness of an insecticide (dieldrin) 
in killing Anopheles farauti mosquitoes. The theory was that resistance to dieldrin was due to a single dominant gene, and that in an appropriately selected sample of the mosquitoes, there should be 50\% susceptibility to insecticide. The hypothesis is 
\begin{equation*}
    H_0:p=\frac{1}{2}\quad \text{versus} \quad H_1:p\neq \frac{1}{2},
\end{equation*}
where $p$ is the probability of susceptibility. The results of such experiment is given in \cite{Osborn}. The sample contains 465 mosquitoes where 264 of them died on being exposed to the insecticide. We can perform this test with test statistic $\text{DDT}_{\beta,\gamma}(\nu_n)$ in Eq.~\ref{eq 38}. Here $\beta$ and $\gamma$ are chosen to be 0.3 and 0.05 respectively. The support of the distribution is $\mathcal{X}=\{0,1\}$, where the digit 1 stands for the death of a mosquito. From here it is evident that $\nu_n(1)=264/465$. The null hypothesis is rejected if the value of the test statistic is large.~In this case the asymptotic null distribution of the test statistic turns out to be $0.774\chi^2_1$. Under the observed data the value of the test statistic turns out to be approximately 6.62.~The 95\% quantile of the aforementioned scaled chi-squared distribution is 2.97.~So, under 5\% level of significance the null hypothesis is rejected. 
{
 \section{Summary}
 In this paper, we have developed a large class of density based divergences which includes the density power divergence family as a special limiting case. The key philosophy of stronger downweighting effect to construct the new family has been discussed. For application purposes, the family gives the data analyst a larger number of choices of possible divergences for inference purposes. We have shown several asymptotic and distributional properties of the proposed estimator. We have also shown that judicial choice of the tuning parameters leads to highly robust and efficient estimators which can often dominate the MDPDE. Though one of the parameters has a smaller effect on the robustness we have shown that both of them play an important role in the context of finite sample efficiency. We have also presented a possible data driven algorithm to obtain the ``optimal'' estimator in a given data set. We have also considered several hypothesis testing strategies for parameteric models which may serve as robust alternatives to the classical likelihood ratio and other likelihood based tests. }
\begin{remark}
Like the MDPDE, the procedures described in this paper avoid the nonparametric density estimation and associated complications specific to classical minimum distance estimation. Another approach of this type can be found in \cite{toma2011dual}.
\end{remark}
\begin{remark}
In creating the test statistics for parametric hypothesis tesing using the
 L$\phi$DPD, we have restricted ourselves to the case where the same set of tuning parameters have been used for estimation as well as the construction of the subsequent divergences. In practice, one could allow them to vary; see, for example, \cite{basu2013testing}. In the present context, while this is possible, we do not explore this issue as we feel that there are enough tuning parameters involved already, and there are no demonstrated results indicating that such differential choices will necessarily produce improved tests.   
\end{remark}

\begin{remark}
In this paper, most of our illustrations have been with respect to the continuous model. Theoretically, however, there is nothing preventing its successful use in discrete models.~All the necessary theories work out satisfactorily in this case. 
\end{remark}
\section{Proof of Theorems}
Proofs of Theorem 2, 5, 6 and 7 are skipped as they can be reproduced along the existing proofs in \cite{Basu}, \cite{ghosh2013robust} and \cite{ghosh2015robustness}. 
\subsection*{Proof of Theorem \ref{thm3} :}
\begin{proof}

(a) From (P4) we know that $f_\theta(x)\phi(f_\theta(x),\gamma)$ is continuous for $\gamma\in(0,1]$ and $\lim_{\gamma\rightarrow 0^+} f_\theta(x)\phi(f_\theta(x),\gamma)=1$. By applying dominated convergence theorem (DCT) on Eq.~\ref{15} at $\gamma\rightarrow 0^+$ we get
$$\frac{1}{n}\sum_{i=1}^{n} u_\theta(X_i) f_\theta^{\beta}(X_i)- \int u_\theta(x) f_\theta^{1+\beta}(x)dx=0$$
which is the unbiased estimating equation for DPD with tuning parameter $\beta$. Hence the result follows.  \\
\\
(b) As $u_\theta(x) f_\theta(x)^{1+\beta}$, $u_{\theta}(x) u_{\theta}(x)^T f_{\theta}(x)^{1+2\beta}$, $u_{\theta}(x) u_{\theta}(x)^T f_{\theta}(x)^{1+\beta}$ are integrable and $f_\theta(x)\phi(f_{\theta}(x),\gamma)$ is bounded, by DCT on Eq.~\ref{16} at $\gamma\rightarrow 0^{+}$ we get
$$J_\beta=\int u_{\theta} u_{\theta}^T f_{\theta}^{1+\beta},\quad K_\beta=\int u_{\theta} u_{\theta}^T f_{\theta}^{1+2\beta}-\zeta_\beta\zeta_\beta^T,$$ $$\zeta_\beta=\int u_{\theta}f_{\theta}^{1+\beta},$$
i.e., $\lim_{\gamma\rightarrow 0^{+}}J_\phi^{-1}K_\phi J_\phi^{-1}=J_\beta^{-1}K_\beta J_\beta^{-1}$. We already know from the assumptions that $J_\beta^{-1}K_\beta J_\beta^{-1}\prec J_\alpha^{-1}K_\alpha J_\alpha^{-1}$, i.e., $(J_\alpha^{-1}K_\alpha J_\alpha^{-1}-J_\beta^{-1}K_\beta J_\beta^{-1})$ is positive definite, where $J_\alpha$ and $K_\alpha$ are defined in the same fashion as $J_\beta$ and $K_\beta$ respectively. The inequality of the asymptotic variances is used here in the sense that AE of DPD with parameter $\beta$ is greater than that of AE of DPD with parameter $\alpha$. So there exists a $\gamma=\gamma_{(\alpha,\beta)}$ such that $J_\phi^{-1}K_\phi J_\phi^{-1}\prec J_\alpha^{-1}K_\alpha J_\alpha^{-1}$. Hence the result follows.
\end{proof} 
\subsection*{Proof of Theorem \ref{thm4} :}
\begin{proof}
First let us assume that breakdown occurs at the model so that there exists sequence $K_n$ of model densities such that $|\theta_n|$ as $n\rightarrow\infty$. Now, consider 
\begin{equation}
D(h_{\epsilon,n},f_{\theta_n})=\int_{A_n} d(h_{\epsilon,n},f_{\theta_n})+\int_{A_n^c}d(h_{\epsilon,n},f_{\theta_n}),
\end{equation}
where $A_n=\{x:g(x)>\max\{k_n(x),f_{\theta_n}(x)\}\}$. Now since $g$ belongs to the model family $\mathcal{F}$, from (BP1) it follows that $\int_{A_n}k_n(x)\rightarrow 0$ and from (BP2) we get $\int_{A_n} f_{\theta_n}\rightarrow 0$, thus under $k_n$ and $f_{\theta_n}$, the set $A_n$ converges to a set of zero probability as $n\rightarrow \infty$.~Thus, on $A_n$, $d(h_{\epsilon,n})\rightarrow d((1-\epsilon)g,0)$ as $n\rightarrow\infty$ and so by DCT $|\int_{A_n}d(h_{\epsilon,n},f_{\theta_n})-\int_{A_n} d((1-\epsilon)g,0)|\rightarrow 0$. Using (BP1), (BP2) and the above result, we have $\int_{A_n}d(h_{\epsilon,n},f_{\theta_n})\rightarrow M_{f,(1-\epsilon)}^{(1)}$. Next, by (BP1) and (BP2), $\int_{A_n^c}g\rightarrow 0$ as $n\rightarrow \infty$, so under $g$, the set $A_n^c$ converges to a set of zero probability.~Hence, similarly, we get $|\int_{A_n^c}d(h_{\epsilon,n},f_{\theta_n})-\int_{A_n^c} d(\epsilon k_n,f_{\theta_n})|\rightarrow 0$. Now by (BP3), we have $\int d(\epsilon k_n,f_{\theta_n})\geq\int d(\epsilon f_{\theta_n},f_{\theta_n})=M_{f,\epsilon}^{(1)}-M_{f,(\epsilon-1)}^{(2)}$. Thus combining the equations we get $\liminf_{n\rightarrow\infty}D(h_{\epsilon,n},f_{\theta_n})\geq M_{f,\epsilon}^{(1)}-M_{f,(\epsilon-1)}^{(2)}+ M_{f,(1-\epsilon)}^{(1)}=a_1(\epsilon),$ say.\\

We will have a contradiction to our breakdown assumption if we can show that there exists a constant value $\theta^*$ in the parameter space such that for the same sequence ${k_n}$, $$\limsup_{n\rightarrow\infty}D(h_{\epsilon,f_{\theta_n}},f_{\theta_n})<a_1(\epsilon)$$ as then the sequence $\{\theta_n\}$ above could not minimize $D(h_{\epsilon,f_{\theta_n}},f_{\theta_n})$  for every $n$. We will now show that above equation is true for all $\epsilon<1/2$ under the model when we choose $\theta^*=\theta^g$.~For any fixed $\theta$, let $B_n=\{x:k_n(x)>\max\{g(x),f_\theta(x)\}\}$. Since $g$ belongs to the model $\mathcal{F}$, from (BP1) we get $\int_{B_n}g\rightarrow 0$, $\int_{B_n}f_\theta\rightarrow 0 $ and $\int_{B_n^c}k_n\rightarrow 0$ as $n\rightarrow\infty$. Thus, under $k_n$, the set $B_n^c$ converges to a set of zero probability, while under $g$ and $f_\theta$, the set $B_n$ converges to a set of zero probability. Thus, on $B_n$, $d(h_{\epsilon,n},f_{\theta})\rightarrow d(\epsilon k_n,0)=B(\epsilon k_n)$ as $n\rightarrow\infty $. So by DCT $|\int_{B_n}d(h_{\epsilon,n},f_{\theta})-\int B(\epsilon k_n)|\rightarrow 0$. Similarly we have $|\int_{B_n^c}d(h_{\epsilon,n},f_{\theta})-\int d((1-\epsilon) g,f_\theta)|\rightarrow 0$. Therefore, we have
\begin{equation}
\limsup_{n\rightarrow\infty}D(h_{\epsilon,n},f_{\theta})=\int D((1-\epsilon)g,f_\theta)+\limsup_{n\rightarrow\infty}\int B(\epsilon k_n).
\label{36}
\end{equation}
Taking $\theta=\theta^g$ in Eq.~\ref{36} and then using (BP3) we get $$\limsup_{n\rightarrow\infty}D(h_{\epsilon,n},f_{\theta^g})\leq M_{f,(1-\epsilon)}^{(1)}-M_{f,(-\epsilon)}^{(2)}+M_{f,\epsilon}^{(1)}=a_3(\epsilon),$$ say. Consequently, asymptotically there is no breakdown if for $\epsilon$ level contamination when $a_3(\epsilon)<a_1(\epsilon)$. But, notice $a_1(\epsilon)$ and $a_3(\epsilon)$ are strictly decreasing and increasing functions respectively. To see this for $a_1(\epsilon)$, notice $M_{f,\epsilon}^{(1)}-M_{f,(\epsilon-1)}^{(2)}=D(\epsilon f_\theta,f_\theta)$. As $\epsilon\uparrow 1$ the above expression decreases. $M_{f,(1-\epsilon)}^{(1)}=\int B((1-\epsilon)f)$. From Eq.~\ref{20} we see that $B(\cdot)$ is an increasing function on positive half line. Using this it is evident that $M_{f,(1-\epsilon)}^{(1)}$ decreases as $\epsilon\uparrow 1$. So, $a_1(\epsilon)$ decreases as $\epsilon\uparrow 1$. Similarly it can be shown $a_3(\epsilon)$ is an increasing function of $\epsilon$. But $a_1(1/2)=a_3(1/2)$; thus asymptotically there is no breakdown and $\limsup_{n\rightarrow\infty}|T_{\beta,\gamma}(H_{\epsilon,n})|<\infty$ for $\epsilon<1/2$. Hence the theorem follows.    
\end{proof}

\subsection*{Proof of Theorem \ref{thm7} :}
\begin{proof}
We know the estimating equation for general M-estimators as \begin{equation*}
\begin{split}
&\frac{1}{n}\sum_{i=1}^nu_{\theta}(X_i)B^{\prime\prime}(f_{\theta}(X_i))f_{\theta}(X_i)-\int u_{\theta}(x)B^{\prime\prime}(f_{\theta}(x))f_{\theta}^2(x)dx\\
&= 0
\end{split}
\end{equation*}
or equivalently \begin{equation*}
\begin{split}
    &\frac{1}{n}\sum_{i=1}^n\Big(u_{\theta}(X_i)B^{\prime\prime}(f_{\theta}(X_i))f_{\theta}(X_i)-\int u_{\theta}(x)B^{\prime\prime}(f_{\theta}(x))f_{\theta}^2(x)dx\Big)\\
    &=0.
\end{split}
\end{equation*}

Viewing this as usual score equation, we take \begin{equation*}
\begin{split}
    &U_\theta(X_i)\\
    &=u_{\theta}(X_i)B^{\prime\prime}(f_{\theta}(X_i))f_{\theta}(X_i)-\int u_{\theta}(x)B^{\prime\prime}(f_{\theta}(x))f_{\theta}^2(x)dx.
    \end{split}
\end{equation*}

We have already seen that the statistic  $T_{\beta,\gamma}(\hat{\theta},\theta_0)$ satisfies $$T_{\beta,\gamma}(\hat{\theta},\theta_0)=n(\hat{\theta}-\theta_0)^TA(\theta_0)(\hat{\theta}-\theta_0)+o_p(1).$$
Note that $\frac{1}{n}\sum_{i=1}^nU_\theta(X_i)=0$ is solved for $\theta=\hat{\theta}$. By Taylor series expansion,
\begin{equation*}
\begin{split}
&\sqrt{n}\frac{1}{n}\sum_{i=1}^nU_{\hat{\theta}}(X_i)\\
&=\sqrt{n}\frac{1}{n}\sum_{i=1}^nU_{\theta_0}(X_i)+\sqrt{n}(\hat{\theta}-\theta_0)\frac{1}{n}\sum_{i=1}^n\frac{\partial}{\partial \theta}U_\theta(X_i)|_{\theta=\theta_0}\\
&+\sqrt{n}(\hat{\theta}-\theta_0)^2\frac{1}{n}\sum_{i=1}^n\frac{\partial^2}{\partial \theta^2}U_\theta(X_i)|_{\theta=\theta^{\prime}}
\end{split}
\end{equation*}
for some $\theta^\prime$ in between $\theta_0$ and $\hat{\theta}$. So we have
\begin{equation*}
\begin{split}
\sqrt{n}\frac{1}{n}\sum_{i=1}^nU_{\hat{\theta}}(X_i) &=\sqrt{n}\bar{U}(\theta_0)\\ &+\sqrt{n}(\hat{\theta}-\theta_0)\frac{1}{n}\sum_{i=1}^n\frac{\partial}{\partial \theta}U_\theta(X_i)|_{\theta=\theta_0}+o_p(1).
\end{split}
\end{equation*}
And hence 
\begin{equation*}
\sqrt{n}\bar{U}(\theta_0)=-\sqrt{n}(\hat{\theta}-\theta_0)\Big[\frac{1}{n}\sum_{i=1}^n\frac{\partial}{\partial \theta}U_\theta(X_i)|_{\theta={\theta_0}}\Big]+o_p(1).
\end{equation*}
Note $\frac{1}{n}\sum_{i=1}^n\frac{\partial}{\partial \theta}U_\theta(X_i)|_{\theta=\theta_0}\to E_{\theta_0}\frac{\partial}{\partial \theta}U_\theta(X_1)|_{\theta=\theta_0}=-J_B(\theta_0)$ as $n\to \infty$.
Hence $$\sqrt{n}\bar{U}(\theta_0)=\sqrt{n}(\hat{\theta}-\theta_0)J_B(\theta_0)+o_p(1).$$
So,
\begin{equation*}
\begin{split}
&n(\hat{\theta}-\theta_0)^TA(\theta_0)(\hat{\theta}-\theta_0)\\ &=n\bar{U}(\theta_0)J_B^{-1}(\theta_0)A(\theta_0)J_B^{-1}(\theta_0)\bar{U}(\theta_0)+o_p(1).
\end{split}
\end{equation*}
This completes the proof.

\end{proof}
\subsection*{Proof of Theorem \ref{thm 8}:}
\begin{proof}
A Taylor expansion of Eq.~\ref{eq 38} around $\hat{\theta}$ gives
\begin{equation}
    \begin{split}
        &\text{DDT}_{\beta,\gamma}(\nu_n)\\ &=2n\left[d_{\beta,\gamma}(\nu_n,f_{\hat{\theta}_0})-d_{\beta,\gamma}(\nu_n,f_{\hat{\theta}})\right]\\
        &=2n\Big[\sum_{j}(\hat{\theta}_{0j}-\hat{\theta}_j)\nabla_j d_{\beta,\gamma}(\nu_n,f_\theta)|_{\theta=\hat{\theta}}\\
        & +\frac{1}{2}\sum_{j,k}(\hat{\theta}_{0j}-\hat{\theta}_j)(\hat{\theta}_{0k}-\hat{\theta}_k)\nabla_{jk}d_{\beta,\gamma}(\nu_n,f_\theta)|_{\theta=\theta^*}\Big]
    \end{split}
    \label{eq 55}
\end{equation}
where the subscripts denote the indicated components of the vector. Also $\theta^*$ lies in the line segment joining $\hat{\theta}_0$ and $\hat{\theta}$. By definition, $\nabla_j d_{\beta,\gamma}(\nu_n,f_\theta)|_{\theta=\hat{\theta}}=0$. Hence, the Eq.~\ref{eq 55} reduces to
\begin{equation}
    \begin{split}
        &\text{DDT}_{\beta,\gamma}(\nu_n)\\&=n(\hat{\theta}_0-\hat{\theta})^T A(\theta_0) (\hat{\theta}_0-\hat{\theta})\\
        &+ n(\hat{\theta}_0-\hat{\theta})^T[\nabla_2 d_{\beta,\gamma}(\nu_n,f_{\theta^*})-A(\theta_0)](\hat{\theta}_0-\hat{\theta})
    \end{split}
    \label{eq 56}
\end{equation}
We will show that under the null $\nabla_2 d_{\beta,\gamma}(\nu_n,f_{\theta^*})\to A(\theta_0)$ as $n\to \infty$. By another Taylor expansion around the true value $\theta_0$, we get for some $\theta^{**}$ between $\theta_0$ and $\theta^*$,
\begin{equation}
\begin{split}
&\nabla_{jk}d_{\beta,\gamma}(\nu_n,f_{\theta^*})\\&= \nabla_{jk}d_{\beta,\gamma}(\nu_n,f_{\theta_0})+\sum_l (\theta_l^*-\theta_{0l})\nabla_{jkl}d_{\beta,\gamma}(\nu_n,f_{\theta^{**}}).
\end{split}
\label{eq 57}
\end{equation}
Under the assumptions (C1)-(C4) it can be easily shown that $\nabla_2 d_{\beta,\gamma}(\nu_n,f_{\theta_0})\to A(\theta_0)$ as $n\to \infty$ and $\nabla_{jkl}d_{\beta,\gamma}(\nu_n,f_{\theta^{**}})=O_P(1)$.~By a simple application of delta theorem on Eq.~\ref{eq 39} it can be shown $\sqrt{n}(\hat{\theta}_0-\theta_0)=O_P(1)$ under the null hypothesis.~Eq.~\ref{eq 40} yields that $\sqrt{n}(\hat{\theta}-\theta_0)=O_P(1)$.~Hence we have $\theta^*=\theta_0+o_P(1)$. As a result the Eq.~\ref{eq 57} reduces to $\nabla_{2}d_{\beta,\gamma}(\nu_n,f_{\theta^*})=A(\theta_0)+o_P(1)$. So, the Eq.~\ref{eq 56} becomes
\begin{equation}
  \text{DDT}_{\beta,\gamma}(\nu_n)=n(\hat{\theta}_0-\hat{\theta})^T A(\theta_0) (\hat{\theta}_0-\hat{\theta})+o_P(1)
  \label{eq 58}
\end{equation}
To obtain the  asymptotic null distribution of $ \text{DDT}_{\beta,\gamma}(\nu_n)$ it is enough to obtain the asymptotic null distribution of $\sqrt{n}(\hat{\theta}_0-\hat{\theta})$. Again from Eq.~\ref{eq 40} and by simple application delta theorem on Eq.~\ref{eq 39} it is easy to show that
$$
\sqrt{n}(\hat{\theta}_0-\hat{\theta})\overset{w}{\to} N(0,\Sigma_{B,b}(\theta_0,\xi_0)),
$$
where $\xi_0$ is the true value of the parameter under $\xi$ formulation of the model. Hence the result follows. If $\Theta_0=\{\theta_0\}$, then Eq.~\ref{eq 58} reduces to 
$$\text{DDT}_{\beta,\gamma}(\nu_n)=n(\theta_0-\hat{\theta})^T A(\theta_0) (\theta_0-\hat{\theta})+o_P(1).$$
We also know 
$$T_{\beta,\gamma}(\hat{\theta},\theta_0)=n(\theta_0-\hat{\theta})^T A(\theta_0) (\theta_0-\hat{\theta})+o_P(1)$$ 
under the null hypothesis.~Hence the asymptotic null distribution of both the statistics are same. This completes the proof.
\end{proof}

\section{Acknowledgements}
The authors gratefully acknowledge the comments of three anonymous referees which led to an improved version of the manuscript. 

\bibliographystyle{alpha} 
\nocite{*}
\bibliography{main}

\end{document}